\documentclass[a4paper,11pt]{amsart}

\usepackage{amssymb,amsthm,amsmath,amscd,amsfonts,bbm,mathrsfs}
\usepackage{stmaryrd}
\usepackage{hyperref}
\usepackage{color}
\hypersetup{
  bookmarksnumbered=true,
  colorlinks=false,
  bookmarksopen=true,
  bookmarksdepth=2,
  pageanchor=true}

\usepackage[cmtip,all]{xy}

\theoremstyle{plain}
\newtheorem{LeABC}{Lemma}
\newtheorem{ThmABC}[LeABC]{Theorem}
\newtheorem{Lemma}[equation]{Lemma}
\newtheorem{Thm}[equation]{Theorem}
\newtheorem{Prop}[equation]{Proposition}
\newtheorem{Cor}[equation]{Corollary}
\theoremstyle{definition}
\newtheorem{Defn}[equation]{Definition}
\newtheorem*{ToDo}{\color{blue}ToDo}
\theoremstyle{remark}
\newtheorem{Remark}[equation]{Remark}
\newtheorem*{Remark*}{Remark}
\newtheorem{Obs}[equation]{Oberservation}
\newtheorem{Notation}[equation]{Notation}
\newtheorem{Terminology}[equation]{Terminology}
\newtheorem*{Notation*}{Notation}
\newtheorem{Ass}[equation]{Assumption}
\newtheorem{Example}[equation]{Example}

\newcommand{\bLe}{\begin{Lemma}}
\newcommand{\eLe}{\end{Lemma}}
\newcommand{\bTh}{\begin{Thm}}
\newcommand{\eTh}{\end{Thm}}
\newcommand{\bPr}{\begin{Prop}}
\newcommand{\ePr}{\end{Prop}}
\newcommand{\bCo}{\begin{Cor}}
\newcommand{\eCo}{\end{Cor}}
\newcommand{\bDe}{\begin{Defn}}
\newcommand{\eDe}{\end{Defn}}
\newcommand{\bAs}{\begin{Ass}}
\newcommand{\eAs}{\end{Ass}}
\newcommand{\bReX}{\begin{Remark*}}
\newcommand{\eReX}{\end{Remark*}}
\newcommand{\bRe}{\begin{Remark}}
\newcommand{\eRe}{\end{Remark}}
\newcommand{\bOb}{\begin{Obs}}
\newcommand{\eOb}{\end{Obs}}
\newcommand{\bEx}{\begin{Example}}
\newcommand{\eEx}{\end{Example}}
\newcommand{\bNo}{\begin{Notation}}
\newcommand{\eNo}{\end{Notation}}
\newcommand{\bTe}{\begin{Terminology}}
\newcommand{\eTe}{\end{Terminology}}
\newcommand{\bNoX}{\begin{Notation*}}
\newcommand{\eNoX}{\end{Notation*}}
\newcommand{\bToDo}{\begin{ToDo}\color{blue}}
\newcommand{\eToDo}{\end{ToDo}}
\newcommand{\btodo}{\bgroup\color{blue}}
\newcommand{\etodo}{\egroup}
\newcommand{\bproof}{\begin{proof}}
\newcommand{\eproof}{\end{proof}}
\newcommand{\beqn}{\begin{equation}}
\newcommand{\eeqn}{\end{equation}}

\renewcommand{\u}{\underline}

\numberwithin{equation}{section}

\usepackage{relsize}
\usepackage[bbgreekl]{mathbbol}

\DeclareSymbolFontAlphabet{\mathbb}{AMSb}
\DeclareSymbolFontAlphabet{\mathbbl}{bbold}

\newcommand{\tcg}[2]{\mathbin{{}_{#1}{\sslash}{}_{#2}}}

\newcommand{\FFF}{\mathcal{F}}

\newcommand{\III}{\mathcal{I}}

\newcommand{\NNN}{\mathcal{N}}

\newcommand{\PPP}{\mathcal{P}}

\newcommand{\SSS}{\mathcal{S}}

\newcommand{\XXX}{\mathcal{X}}
\newcommand{\YYY}{\mathcal{Y}}
\newcommand{\ZZZ}{\mathcal{Z}}

\newcommand{\Fb}{\mathfrak{b}}
\newcommand{\Fg}{\mathfrak{g}}

\newcommand{\Fu}{\mathfrak{u}}

\newcommand{\FD}{\mathfrak{D}}

\newcommand{\FZ}{\mathfrak{Z}}

\newcommand{\FF}{{\mathbb{F}}}
\newcommand{\GG}{{\mathbb{G}}}

\newcommand{\ZZ}{{\mathbb{Z}}}

\DeclareMathOperator{\colim}{colim}

\DeclareMathOperator{\conj}{conj}
\DeclareMathOperator{\can}{can}
\DeclareMathOperator{\crys}{cris}
\DeclareMathOperator{\eff}{eff}

\DeclareMathOperator{\et}{\acute et}
\DeclareMathOperator{\fppf}{fppf}
\DeclareMathOperator{\gr}{gr}
\DeclareMathOperator{\id}{id}
\DeclareMathOperator{\indsyn}{ind-syn}

\DeclareMathOperator{\nil}{nil}

\DeclareMathOperator{\perf}{perf}

\DeclareMathOperator{\qrsp}{qrsp}
\DeclareMathOperator{\qsyn}{qsyn}
\DeclareMathOperator{\rad}{rad}
\DeclareMathOperator{\real}{real}

\DeclareMathOperator{\res}{res}

\DeclareMathOperator{\syn}{syn}

\DeclareMathOperator{\Ab}{Ab}
\DeclareMathOperator{\Ad}{Ad}
\DeclareMathOperator{\Aff}{Aff}
\DeclareMathOperator{\Aut}{Aut}
\DeclareMathOperator{\uAut}{\underline{Aut}}

\DeclareMathOperator{\BT}{BT}
\DeclareMathOperator{\bBT}{\overline{BT}}

\DeclareMathOperator{\Disp}{Disp}
\DeclareMathOperator{\tDisp}{\widetilde{Disp}}

\DeclareMathOperator{\Fil}{Fil}
\DeclareMathOperator{\Fun}{Fun}

\DeclareMathOperator{\GL}{GL}
\DeclareMathOperator{\Hom}{Hom}

\DeclareMathOperator{\LG}{LG}
\DeclareMathOperator{\Lau}{Lau}
\DeclareMathOperator{\Lie}{Lie}

\DeclareMathOperator{\Mod}{Mod}
\DeclareMathOperator{\Nil}{Nil}
\DeclareMathOperator{\tNil}{\widetilde{Nil}}

\DeclareMathOperator{\Spf}{Spf}

\DeclareMathOperator{\Spec}{Spec}

\DeclareMathOperator{\pWin}{pWin}
\DeclareMathOperator{\Win}{Win}

\newcommand{\ttDisp}{{\tDisp}{}}
\newcommand{\ttNil}{{\tNil}{}}
\newcommand{\bbBT}{{\bBT}{}}
\newcommand{\sN}{{}^sN}
\newcommand{\sI}{{}^sI}
\newcommand{\sW}{{}^sW}
\newcommand{\srho}{{}^s\!\rho}
\newcommand{\svarrho}{{}^s\!\varrho}

\begin{document}

\title{The Shimurian BT stack is a gerbe over truncated displays}
\author{Eike Lau}
\date{\today}
\address{Eike Lau, Fakult\"{a}t f\"{u}r Mathematik,
Universit\"{a}t Bielefeld, D-33501 Bielefeld}

\begin{abstract}
We show that the mod $p$ fiber of the Shimurian stack $\BT_n^{G,\mu}$ constructed by Gardner--Madapusi is a gerbe over the corresponding stack of truncated displays. This confirms a conjecture of Drinfeld.
\end{abstract}

\maketitle

\setcounter{tocdepth}{1}
\tableofcontents

\section{Introduction}

\subsection{Main result}

The fibered category of $n$-truncated Barsotti--Tate groups is an algebraic stack $\BT_n$, which is smooth of relative dimension zero over $\Spec \ZZ$ by a theorem of Grothendieck proved in \cite{Illusie-BT}. The reduction $\bbBT_n=\BT_n\times\Spec\FF_p$ admits a morphism $\phi_n:\bbBT_n\to\Disp_n$ to the stack of $n$-truncated displays, and $\phi_n$ is a gerbe banded by a certain finite locally free group scheme, by \cite{Lau:Smoothness}. This group scheme was determined explicitly in \cite{Drinfeld:The-Lau}.

According to a series of conjectures developed in \cite{Drinfeld:On-Shimurian, Drinfeld:The-Lau}, this picture has the following Shimurian analogue. Given a smooth affine group scheme $G$ over $\ZZ_p$ and a $1$-bounded cocharacter $\mu$ of $G$ defined over $W(k)$, for a finite extension $k$ of $\FF_p$, there exists a smooth $p$-adic formal algebraic stack $\BT_n^{G,\mu}$ over $\Spf W(k)$ together with a morphism 
\beqn
\label{Eq:phinGmu}
\phi_n^{G,\mu}:\bbBT_n^{G,\mu}\to\Disp_n^{G,\mu}
\eeqn
from the special fiber to the stack of $n$-truncated $(G,\mu)$-displays, the morphism $\phi_n^{G,\mu}$ is a gerbe banded by an explicit finite locally free group scheme\footnote{\label{Fn:FD-Lau}The group scheme $\FD_n^{G,\mu}$ is denoted $\Lau_n^{G,\mu}$ in \cite{Drinfeld:The-Lau}.} $\FD_n^{G,\mu}$ over $\Disp_n^{G,\mu}$, and setting $G=\GL_h$ for varying $h$ and minuscule cocharacter $\mu$ gives back the $\BT_n$ case. 

The following parts of this picture are known. 
The stack $\BT_n^{G,\mu}$ is defined in \cite{Drinfeld:On-Shimurian,Gardner-Madapusi}\footnote{The definition in \cite{Drinfeld:On-Shimurian} assumes that $k=\FF_p$; the general definition is given in \cite{Gardner-Madapusi}.} 
as a derived $p$-adic formal stack using the prismatization of \cite{Bhatt-Lurie:The-prismatization, Drinfeld:Prismatization}. It is a classical smooth $p$-adic formal algebraic stack by \cite[Theorem D]{Gardner-Madapusi}, which proves the algebraicity conjecture \cite[Conjecture C.3.1]{Drinfeld:On-Shimurian}. The $p$-adic formal algebraic stack $\BT_n^{h,d}$ over $\ZZ_p$ of $n$-truncated Barsotti-Tate groups of height $h$ and tangent dimension $d$ is isomorphic to $\BT_n^{\GL_h,\mu_d}$, where $\mu_d$ is the minuscule cocharacter with determinant $t\mapsto t^d$, by \cite[Theorem 11.2.7]{Gardner-Madapusi}. The morphism $\phi_n^{G,\mu}$ is provided by \cite[Remark 9.1.3]{Gardner-Madapusi}.

In this article, we will add the following missing piece.

\begin{ThmABC}
\label{Th:Main}
The morphism $\phi_n^{G,\mu}$ is a gerbe banded by the group scheme $\FD_n^{G,\mu}$.
\end{ThmABC}

This confirms the gerbe conjecture \cite[Conjecture C.5.3]{Drinfeld:The-Lau}; see footnote \ref{Fn:FD-Lau}.
The case $n=1$ is known by the proof of \cite[Theorem 9.3.2]{Gardner-Madapusi}; cf.\ \cite[C.5.4, Remark (i)]{Drinfeld:The-Lau}.

\subsection{Reduction to an elementary version}

Since the source and target of $\phi_n^{G,\mu}$ are smooth algebraic stacks, it suffices to look at the functors $\phi_n^{G,\mu}(R)$ for quasiregular semiperfect (qrsp) $k$-algebras $R$. These functors have the following classical description. The torsion-free ring $A_{\crys}(R)$ with the Nygaard filtration gives a frame $\u A_{\crys}(R)$ over $W(k)$ in the sense of \cite{Lau:Higher}. Let $\u W(R)$ be the Witt frame. There is a natural homomorphism $\varrho:\u A_{\crys}(R)\to \u W(R)$, which induces a homomorphism of the reductions modulo $p^n$, denoted $\varrho_n:\u A_n(R)\to\u {W\!}_n(R)$, and $\phi_n^{G,\mu}(R)$ is the associated functor of groupoids of $(G,\mu)$-windows\footnote{\label{Fn:term-display-intro}These objects were called $G$-displays of type $\mu$ in \cite{Lau:Higher}. 
} as defined in loc.cit.

Let us recall the format of the definition of $(G,\mu)$-windows. For any frame $\u C$ over $W(k)$, one defines groups and group homomorphisms
$\sigma,\tau:G(C^\oplus)_\mu\rightrightarrows G(C_0)$
and a right action of the first group on the second group, viewed as a set, given by $x*y=\tau(y)^{-1}x\sigma(y)$. We denote by $\pWin^{G,\mu}(\u C)$ the quotient groupoid for this action; this is a presheaf version of the groupoid of $(G,\mu)$-windows over $\u C$.  Now the functor $\phi_n^{G,\mu}(R)$ for a qrsp $k$-algebra $R$ arises by etale stackification from the functor of groupoids 
\beqn
\label{Eq:varrhonGmu-intro}
\varrho_n^{G,\mu}(R):\pWin^{G,\mu}(\u A_n(R))\to\pWin^{G,\mu}(\u {W\!}_n(R))
\eeqn
induced by $\varrho_n$. The following is the central technical result of this article.

\begin{ThmABC}
\label{Th:Main-pre}
The functor $\varrho_n^{G,\mu}$ is surjective on objects and full, with inertia given by $\FD_n^{G,\mu}$.
\end{ThmABC}

Some comments on the proof of Theorem \ref{Th:Main-pre} are given below.
Theorem \ref{Th:Main} is a formal consequence of Theorem \ref{Th:Main-pre} since a morphism of smooth algebraic stacks over $k$ is a gerbe banded by a given syntomic group scheme iff this property holds for the restriction to qrsp $k$-algebras. By an elaboration of such arguments, Theorem \ref{Th:Main-pre} implies directly that the etale stackification of $\varrho_n^{G,\mu}(R)$ for varying qrsp $k$-algebras $R$ extends to a unique morphism of algebraic stacks $\ttDisp_n^{G,\mu}\to\Disp_n^{G,\mu}$ over $k$. Since the derived stack $\bbBT_n^{G,\mu}$ is a classical smooth algebraic stack by the results of \cite{Gardner-Madapusi}, it follows that the stack $\ttDisp_n^{G,\mu}$ is isomorphic to $\bbBT_n^{G,\mu}$ because they agree on qrsp $k$-algebras.

\subsection{Drinfeld's group scheme}

Let us briefly recall the group scheme $\FD_n^{G,\mu}$ over $\Disp_n^{G,\mu}$ defined in \cite{Drinfeld:The-Lau} that appears in Theorems \ref{Th:Main} and \ref{Th:Main-pre}.

A central object of Zink's theory of displays in \cite{Zink:The-Display} is a functor $\FZ:\Disp^{\{0,1\}}\to\LG$ where $\Disp^{\{0,1\}}$ is the category of classical displays over $p$-nilpotent rings, viewed here as higher displays of type $\{0,1\}$, and $\LG$ is the category of commutative formal Lie groups. Originally the functor was denoted $BT$.
This construction allows the following extension. Let $\Disp_n^{\eff}$ denote the category of effective $n$-truncated displays over $\FF_p$-algebras; effective means type in $[0,m]$ for some $m$. Let $\LG_n$ be the category of commutative $n$-smooth group schemes over $\FF_p$-algebras as defined in \cite{Drinfeld:The-Lau}; these are finite group schemes which arise locally as the $F^n$-kernel of a commutative formal Lie group. Then there is a functor
\[
\FZ_n:\Disp_n^{\eff}\to\LG_n
\]
defined by the same formula as $\FZ$. The extension of $\FZ$ to truncated objects was observed in \cite{Lau-Zink}, and the extension to higher effective (truncated) displays appears in \cite{Drinfeld:The-Lau}.

The adjoint representation of $G$ gives a functor $\Lie^{G,\mu}:\Disp_n^{G,\mu}\to\Disp_n$. If $\mu$ is $1$-bounded, the twisted functor $\Lie^{G,\mu}(-1)$ has image in $\Disp_n^{\eff}$, and we have $\FD_n^{G,\mu}=\FZ_n\circ\Lie^{G,\mu}(-1)$ as functors $\Disp_n^{G,\mu}\to\LG_n$.

\subsection{Sheared Witt vectors}

For a semiperfect $\FF_p$-algebra $R$ we will need modifications of the ring of Witt vectors $W(R)$ defined by
\[
\sW^{(n)}(R)=W(R)/\hat W(R[F^n])
\] 
and $\sW(R)=\lim_n\sW^{(n)}(R)$ with transition maps given by the Witt vector Frobenius. Equivalently, if $R=R^\flat/J$ where $R^\flat$ is the limit perfection, then $\sW(R)=W(R^\flat)/\hat W(J)$; here $\hat W(J)$ is the set of Witt vectors with coefficients in $J$ that converge to zero for the limit topology of $R^\flat$. 
If $R$ is a qrsp $\FF_p$-algebra, the ring $\sW(R)$ is related to $A_{\crys}(R)$ as follows. The ideal $N=\ker(A_{\crys}(R)\to W(R))$ carries the $\varphi$-linear endomorphism $\dot\varphi=p^{-1}\varphi$. Let $N^{\nil}$ be the set of elements of $N$ on which $\dot\varphi$ is $p$-adically nilpotent.

\begin{ThmABC}
\label{Th:Acris-sW}
If $R$ is a qrsp $\FF_p$-algebra, there is an isomorphism $\sW(R)\cong A_{\crys}(R)/N^{\nil}$.
\end{ThmABC}

The ring $\sW(R)$ can be defined for every $p$-nilpotent ring $R$ and gives an fpqc sheaf $\sW$ of rings studied in \cite{Dr25:Ring-stacks-conjecturally}, and in \cite{BMVZ} as a basic ingredient of a theory of sheared prismatization developed by the authors of loc.\,cit.\ together with A.~Kanaev. In the first draft of this article, $\sW(R)$ was denoted $\check W(R)$. We adopt the terminology and notation of loc.\,cit.

\subsection{On the proof of Theorem \ref{Th:Main-pre}}
\label{Se:Intro.inertia}

Let $R$ be a qrsp $k$-algebra.
We will use auxiliary frames whose underlying rings fit into a commutative diagram of surjective homomorphisms
\[
\xymatrix@M+0.2em{
A_{\crys}(R) \ar[r] & B_n(R) \ar[r]^-{\svarrho_n} \ar[d] & \sW^{(n)}(R) \ar[r] & W(R) \ar[d] \\
&  A_n(R) \ar[rr]^{\varrho_n} && W_n(R)
}
\]
with $B_n(R)=A_{\crys}(R)/p^n N$. One deduces an exact sequence of $B_n(R)$-modules
\[
0\to\hat W(R[F^n])\to\ker(\svarrho_n)\to\ker(\varrho_n)\to 0
\]
where $\ker(\svarrho_n)$ can be identified with the set of $\dot\varphi$-nilpotent elements in $N/p^n$ using Theorem \ref{Th:Acris-sW}. With a frame version of this sequence, the classical deformation theory of $(G,\mu)$-windows yields that the functor $\varrho_n^{G,\mu}(R)$ of \eqref{Eq:varrhonGmu-intro} is full, and its inertia at a given element $g\in G(A_n(R))$ with image $\bar g\in G(W_n(R))$ can be identified with the cokernel of a homomorphism of abelian groups
\[
(\Lie(G)\otimes\hat W(R)^\oplus)_\mu\xrightarrow{\;\Ad(\bar g)\otimes\varphi-{\id}\otimes\tau\;}\Lie(G)\otimes\hat W(R)
\]
that arises from a Lie version of groupoid of $(G,\mu)$-windows. By the definition of $\FD_n^{G,\mu}$, this cokernel is the group of sections of $\FD_n^{G,\mu}$ over the morphism $\bar g:\Spec R\to\Disp_n^{G,\mu}$. 

%
%
%

\subsection{Outline}

Sections \S\ref{Se:frames-modules-ideals}--\ref{Se:Def-Gmu-groupoids} are mainly concerned with notation and recollections from the literature with some simple additions. 
In \S\ref{Se:Witt-crys-frames}, \ref{Se:Witt-crys-Gmu-disp} we introduce the functor \eqref{Eq:varrhonGmu-intro} and state the central Proposition \ref{Pr:Main}, which is a preliminary version of Theorem \ref{Th:Main-pre}. Sheared Witt vectors and Theorem \ref{Th:Acris-sW} are discussed in \S\ref{Se:Sheared}, then Proposition \ref{Pr:Main} is proved in \S\ref{Se:frame-B}, \ref{Se:Proof-Prop}, and Theorem \ref{Th:Main-pre} is deduced from this in \S\ref{Se:Drinfeld-group-scheme} after recollections on the functors $\FZ$ and $\FZ_n$ in \S\ref{Se:Zink-functor}.
The final applications including Theorem \ref{Th:Main} are collected in \S\ref{Se:Geom-appl}, using generalities on the reconstruction of algebraic stacks from quasiregular semiperfect rings in \S\ref{Se:Reconstruction}.

\subsection{Acknowledgements}

The author is grateful to V.~Drinfeld for helpful discussions and generously sharing unpublished notes and drafts. The author thanks M.~Hoff and T.~Zink for valuable conversations.
This work was 
funded by the Deutsche Forschungsgemeinschaft (DFG, German Research
Foundation) -- Project-ID 491392403 -- TRR 358.

\section{Frames and modules}
\label{Se:frames-modules-ideals}

We recall the notion of (higher) frames of \cite{Lau:Higher} and fix some additional terminology. An animated version of frames was recently introduced in \cite{Gardner-Madapusi}, but here we can stay classical.

\subsection{Frames}
\label{Se:Frames.Terminology}

A frame $\u A=(A,\sigma,\tau)$ consists of a $\ZZ$-graded ring $A=\bigoplus_i A_i$ and ring homomorphisms $\sigma,\tau:A\to A_0$ with the following properties. Let $\sigma_n,\tau_n:A_n\to A_0$ be the restrictions of $\sigma$ and $\tau$. Then $\tau_0$ is the identity, $\tau_n$ is bijective for $n<0$, the endomorphism $\sigma_0$ of $A_0$ is a Frobenius lift, $\sigma_{-1}(t)=p$ where $t\in A_{-1}$ is the unique element with $\tau_{-1}(t)=1$, and $p\in\rad(A_0)$. 
Let $R=A_0/tA_1$. We say that $\u A$ is a frame for $R$. 

\bRe
We have $\tau_n(a)=t^na$ for $n\in\ZZ$ and $a\in A_n$; the formula makes sense for all $n$ because $t:A_0\to A_n$ is bijective for $n\ge 0$.
\eRe

\subsubsection*{Torsion and complete frames}
A frame $\u A$ will be called torsion if $A$ is annihilated by a power of $p$, and $\u A$ will be called (derived) $p$-complete if all $A_i$ are (derived) $p$-complete. Note that for any derived $p$-complete ring $B$ the pair $(B,pB)$ is Henselian by \cite[\href{https://stacks.math.columbia.edu/tag/0G3H}{Lemma 0G3H}]{Stacks}, in particular $p\in\rad(B)$. 
The frames of interest in this article will be torsion frames, arising as quotients of $p$-complete frames, but some auxiliary frames will only be derived $p$-complete.

\subsubsection*{$1$-frames}
The following variant of the frames of \cite{Lau:Frames} will be called $1$-frames: A $1$-frame is a quadruple $(A,I,\sigma,\dot\sigma)$ consisting of a ring $A$ with $p\in\rad(A)$, an ideal $I\subseteq A$, a Frobenius lift $\sigma:A\to A$, and a $\sigma$-linear map $\dot\sigma:I\to A$ such that $p\dot\sigma=\sigma|_I$.
A frame $\u A=(A,\sigma,\tau)$ such that $\tau_1:A_1\to A_0$ is injective gives a $1$-frame $(A,I,\sigma,\dot\sigma)$ with $I=\tau(A_1)$ and $\dot\sigma=\sigma_1\circ\tau_1^{-1}$.

\subsection{Modules}

Let $\u A=(A,\sigma,\tau)$ be a frame. 

\bDe 
An $\u A$-module $\u M=(M,F)$ consists of a graded $A$-module $M$ and a homomorphism of $A_0$-modules $F:M^\sigma\to M^\tau$.  An $\u A$-window is an $\u A$-module $(M,F)$ such that $M$ is a finitely generated projective $A$-module and $F$ is bijective.\footnote{$\u A$-modules and $\u A$-windows were called predisplays and displays over $\u A$ in \cite[Definition 3.2.1]{Lau:Higher}. 
} 
\eDe


\subsubsection*{Categories of modules and windows}

$\u A$-modules form an abelian category $\Mod(\u A)$ because the functor $M\mapsto M^\tau=M/(t-1)$ is exact. $\u A$-windows form an exact subcategory $\Win(\u A)$. There is a tensor product of $\u A$-modules defined by $(M,F)\otimes(M',F')=(M\otimes_AM',F\otimes F')$, and the tensor product of windows is a window. 
The shift of a graded $A$-module $M$ by $m\in\ZZ$ is defined by $M(m)_n=M_{m+n}$; this gives a shift of $\u A$-modules and $\u A$-windows defined by $(M,F)(m)=(M(m),F)$, using that $M(m)^\sigma=M^\sigma$ and $M(m)^\tau=M^\tau$.

\subsubsection*{Structure of modules}
 
A graded $A$-module $M$ will be called effective if the multiplication map $t^i:M_0\to M_{-i}$ is bijective for $i\ge 0$.
Each finitely generated projective graded $A$-module $M$ takes the form $M\cong L\otimes_{A_0}A$ where $L$ is a finite projective graded $A_0$-module, thus $L=\bigoplus_i L_i$ and $M_n=\bigoplus_i L_i\otimes_{A_0}A_{n-i}$ in degree $n$. Here $L\otimes_{A_0}A$ is effective iff $L_i=0$ for $i<0$. 

\subsubsection*{Functoriality}

For a frame homomorphism $\u A\to\u A'$ there is a restriction of scalars functor $\Mod(\u A')\to\Mod(\u A)$ since the functor $M\mapsto M^\tau=M/(t-1)$ is compatible with restriction of scalars. The restriction of scalars has a left adjoint $\Mod(\u A)\to\Mod(\u A')$ given by $(M,F)\mapsto(M\otimes_AA',F\otimes\sigma)$, and this functor restricts to a functor $\Win(\u A)\to\Win(\u A')$.

\subsubsection*{Fundamental complex}

For an $\u A$-module $\u M=(M,F)$ and $i\in\ZZ$ we will consider the complex of abelian groups in degrees $\{-1,0\}$
\beqn
\label{Eq:Ci}
\Gamma_i(\u M)=[M_i\xrightarrow{\gamma_F} M^\tau]
\eeqn
where $\gamma_F(y)=F(y^\sigma)-y^\tau$ with $y^\sigma=y\otimes 1$ in $M^\sigma$ and $y^\tau=y\otimes 1$ in $M^\tau$.
If $M$ is effective, then the homomorphism $M_0\to M^\tau$, $x\mapsto x^\tau$ is bijective, thus $F$ corresponds to a $\sigma$-linear map $\phi:M\to M_0$, and 
\beqn
\label{Eq:Ci-eff}
\Gamma_i(\u M)=[M_i\xrightarrow{\;\phi_i-t^i\;} M_0].
\eeqn

\subsection{Ideals}

Let $\u A=(A,\sigma,\tau)$ be a frame. 

\bDe
\label{Def:ideal-frame}
A homogeneous ideal $K\subseteq A$ will be called an ideal of the frame $\u A$ if the graded ring $A/K$ carries an induced frame structure, denoted by $\u A/K=(A/K,\bar\sigma,\tau)$. Here induced means that the projection map $A\to A/K$ is a frame homomorphism $\u A\to\u A/K$. 
\eDe

\bRe
\label{Re:ideal-frame}
A homogeneous ideal $K$ of $A$ is an ideal of the frame $\u A$ iff the ring homomorphisms $\sigma$ and $\tau$ restrict to homomorphisms $K\to K_0$ and the map $t^i:K_0\to K_{-i}$ is bijective for $i\ge 0$. 
\eRe

\bEx
Every frame $\u A$ carries the ideals $p^nA$, which give torsion frames $\u A/p^n$.
\eEx

\bRe
\label{Re:ideal-module}
Each ideal $K$ of $\u A$ gives rise to an $\u A$-module $\u K=(K,F_{\can})$ as a submodule of the tautological $\u A$-module $(A,\id_{A_0})$.
\eRe

\bDe
\label{De:leveled}
An ideal $K$ of $\u A$ will be called leveled if $\tau_1^K:K_1\to K_0$ is bijective. In that case let $\dot\sigma=\dot\sigma_K=\sigma_1\circ(\tau_1^K)^{-1}$; this is a $\sigma_0$-linear endomorphism of $K_0$. A leveled ideal $K$ will be called $\dot\sigma$-nilpotent if $\dot\sigma:K_0/p\to K_0/p$ is locally nilpotent.
\eDe

\section{$(G,\mu)$-groupoids}
\label{Se:Gmu-groupoid}

We fix a pair $(G,\mu)$ where $G$ is a smooth affine group scheme over $\ZZ_p$ and $\mu$ is a cocharacter of $G$ defined over $W(k)$ for a finite extension $k$ of $\FF_p$. 
Let $\u A$ be a frame as in \S\ref{Se:Frames.Terminology} such that $A$ is a graded $W(k)$-algebra. 
Following \cite[\S 5.1]{Lau:Higher}, to this data we associate a groupoid of $(G,\mu)$-windows (called displays in loc.cit.), here in a presheaf version.

\subsection{Preliminaries}

\subsubsection*{Twisted conjugation}

Assume that groups $H$, $H'$ and group homomorphisms $\sigma,\tau:H\to H'$ are given. We define a right action of $H$ on the set $H'$ by 
$g*h=\tau(h)^{-1}g\sigma(h)$
for $g\in H'$ and $h\in H$. The associated action groupoid 
\beqn
\label{Eq:twisted-conj-groupoid}
(H'\tcg{\sigma}{\tau}H)
\eeqn
will be called the twisted conjugation groupoid associated to $(H\rightrightarrows H')$. The set of objects of $(H'\tcg{\sigma}{\tau}H)$ is $H'$, and $h\in H$ gives a morphism $g*h\to g$. The notation $\tcg{\sigma}{\tau}$ is borrowed from \cite[\S5.5]{Gardner-Madapusi}.
If the groups $H$, $H'$ are abelian, then $(H'\tcg{\sigma}{\tau}H)$ is the underlying groupoid of the Picard category associated to the complex of abelian groups $[H\xrightarrow{\sigma-\tau}H']$  in degrees $\{-1,0\}$.

\subsubsection*{Sheaf version}

If $H\rightrightarrows H'$ are sheaves of groups over a site $(\SSS,T)$, then $(H'\tcg{\sigma}{\tau}H)$ is a sheaf of groupoids in the $1$-categorical sense, which gives a prestack over $\SSS$. The associated stack will be denoted by $[H'\tcg{\sigma}{\tau}H]$ or $[H'\tcg{\sigma}{\tau}H]_T$.

\subsubsection*{Actions of $\GG_m$}

Assume that $X=\Spec C$ is an affine scheme over a ring $R$ with a right action of $\GG_m$. The action corresponds to a $\ZZ$-grading $C=\bigoplus_iC_i$ as follows. We let $\lambda\in\GG_m$ act on $f\in C$ by $(\lambda\cdot f)(x)=f(x\cdot \lambda)$. Then $C_i$ is the set of all $f\in C$ such that $\lambda\cdot f=\lambda^if$. 
For a graded $R$-algebra $B$, thus a right action of $\GG_m$ on $\Spec B$, the set of $\GG_m$-equivariant morphisms $\Spec B\to X$ over $R$, or equivalently of homomorphisms of graded $R$-algebras $C\to B$, will be denoted by $X(B)^{\GG_m}$.

\subsection{Definition of the groupoid}
\label{Se:displays-Gmu-groupoid}

Let $(G,\mu)$ and $\u A$ be as above. Let $\GG_m$ act on $G_{W(k)}$ by right conjugation with $\mu$, so $\lambda\in\GG_m$ acts on $g\in G_{W(k)}$ by $g*\lambda=\mu(\lambda)^{-1}g\mu(\lambda)$.
For a $\GG_m$-invariant closed subscheme $X\subseteq G_{W(k)}$ let
\[
X(A)_\mu=X(A)^{\GG_m}.
\]
The homomorphism $W(k)$-algebras $\tau:A\to A_0$ induces a map 
\[
\tau:X(A)_\mu\to X(A_0).
\]
If $X$ is defined over $\ZZ_p$, the ring homomorphism $\sigma:A\to A_0$ induces a map 
\[
\sigma:X(A)_\mu\to X(A_0).
\]
In particular, we have group homomorphisms $\sigma,\tau:G(A)_\mu\to G(A_0)$. The associated twisted conjugation groupoid \eqref{Eq:twisted-conj-groupoid} will be denoted by
\beqn
\label{Eq:Gamma-GmuA}
\pWin^{G,\mu}(\u A)=(G(A_0)\tcg{\sigma}{\tau}G(A)_\mu)
\eeqn
and will be called the $(G,\mu)$-groupoid over $\u A$. A sheaf version of this construction for the etale topology was called the groupoid of $G$-displays of type $\mu$ over $\u A$ in \cite[Definition 5.3.2]{Lau:Higher}.

\subsection{Decompositions}
\label{Se:displays-decompositions}

Let $P^+\subseteq G_{W(k)}$ and $P^-\subseteq G_{W(k)}$ be the maximal closed subschemes on which $\GG_m$ acts with non-negative and non-positive weights. If $G_{W(k)}=\Spec C$ then $P^\pm=\Spec C^\pm$ where $C^+=C/(C_{<0})$ and $C^-=C/(C_{>0})$. 
Then $P^{\pm}$ are $\GG_m$-stable subgroups of $G_{W(k)}$, and $P^+\cap P^-$ is the subgroup of $\GG_m$-invariants. There are unique $\GG_m$-equivariant retractions $P^\pm\to P^+\cap P^-$; let $U^\pm$ be the kernels. Multiplication gives a bijective map
\beqn
\label{Eq:GAmu-decomposition}
P^-(A)_\mu\times U^+(A)_\mu\xrightarrow\sim G(A)_\mu
\eeqn
by \cite[Proposition 6.2.2]{Lau:Higher}, and the homomorphism
\beqn
\label{Eq:Pamu-PA0}
\tau:P^-(A)_\mu\xrightarrow\sim P^-(A_0)
\eeqn
is an isomorphism by \cite[Lemma 6.2.1]{Lau:Higher}. 

\subsubsection*{Lie algebra and weights}

Let $\Fg=\Lie(G)$, a finite free $\ZZ_p$-module, and $\Fu^+=\Lie(U^+)$, a finite free $W(k)$-module. The given right action of $\GG_m$ on $G_{W(k)}$ induces a right action of $\GG_m$ on $\Fg\otimes_{\ZZ_p} W(k)$ and hence a $\ZZ$-grading\footnote{Here the upper numbering is used since this grading is associated to a right action of $\GG_m$, as opposed to the left action of $\GG_m$ on $A$ and on the coordinate ring of $G$.} 
\beqn
\label{Eq:grading-Fg-Wk}
\Fg\otimes_{\ZZ_p} W(k)=\bigoplus_{i\in\ZZ}\Fg^i
\eeqn 
under which $\Fu^+=\bigoplus_{i>0}\Fg^i$. There is a $\GG_m$-equivariant isomorphism $\log_{U^+}:U^+\xrightarrow\sim V(\Fu^+)$ of schemes over $W(k)$ 
which induces the identity on Lie algebras where $V$ means vector group; see \cite[Lemma 6.1.1]{Lau:Higher}. The isomorphism $\log_{U^+}$ induces a bijective map 
\beqn
\label{UAmu-uAi}
U^+(A)_\mu\xrightarrow\sim\prod_{i>0}\Fg^i\otimes_{W(k)}A_i.
\eeqn
Following \cite[Definition 6.3.1]{Lau:Higher} the pair $(G,\mu)$ will be called $1$-bounded if $\Fu^+=\Fg^1$. In that case, $\log_{U^+}$ is unique and is an isomorphism of group schemes by \cite[Lemma 6.3.2]{Lau:Higher}, and \eqref{UAmu-uAi} is an isomorphism of groups
\beqn
\label{UAmu-uA1}
U^+(A)_\mu\xrightarrow\sim\Fu^+\otimes_{W(k)}A_1.
\eeqn

\bEx
Let $\mu:\GG_m\to G=\GL_2$ be $\mu(\lambda)=\left(\begin{smallmatrix}\lambda \\ & 1\end{smallmatrix}\right)$. Then $G(A)_\mu$ is the group of  invertible matrices $\left(\begin{smallmatrix}a&b\\c&d\end{smallmatrix}\right)$ with $a,d\in A_0$, $b\in A_{-1}$, and $c\in A_1$. The subgroup $P^-(A)_\mu$ is defined by $c=0$, and $U^+(A)_\mu$ consists of all $\left(\begin{smallmatrix}1&0\\c&1\end{smallmatrix}\right)$.
\eEx

\subsection{The case \boldmath $G=\GL(V)$}
\label{Se:Groupoids.GL(V)}

Let $V$ be a finite free $\ZZ_p$-module and consider $\GL(V)$ as a group scheme over $\ZZ_p$. Let $\tilde\mu$ be a cocharacter of $\GL(V)$ defined over $W(k)$ and let $V_{\tilde\mu}=\bigoplus_iV_{\tilde\mu,i}$ be the associated graded $W(k)$-module, thus $V_{\tilde\mu}=V_{W(k)}$ as a module, and $\GG_m$ acts on $V_{\tilde\mu,i}$ via $\lambda\mapsto\lambda^{i}$. 
If we view $V_{\tilde\mu,i}$ as a graded module in degree zero, the formula reads $V_{\tilde\mu}=\bigoplus_iV_{\tilde\mu,i}(-i)$.
\bLe
\label{Le:GLn-groupoid}
There is a realisation functor
$\real:\pWin^{\GL(V),\tilde\mu}(\u A)\to\Win(\u A)$
which induces an isomorphism between the source and the groupoid of $\u A$-windows $(M,F)$ with $M=V_{\tilde\mu}\otimes_{W(k)}A$ as a graded $A$-module.
\eLe

\bproof
See \cite[Example 5.3.5]{Lau:Higher}.
Let $N=V_{\tilde\mu}\otimes_{W(k)}A$ as a graded $A$-module.
For $g\in\GL(V)(A_0)$ let $\real(g)=(N,F_g)$ where $F_g$ is given by $g:V_{A_0}\to V_{A_0}$ under the natural identification $N^\tau=V_{A_0}=N^\sigma$.
The group $\GL(V)(A)_{\tilde \mu}$ is the group of graded module automorphisms of $N$, and an element $h$ of this group defines an automorphism of $\u A$-modules $(N,F_{g*h})\to(N,F_g)$ where $g*h=\tau(h)^{-1}g\sigma(h)$; this automorphism is $\real(h)$. 
\eproof

\section{Variations of $(G,\mu)$-groupoids} 

Let $(G,\mu)$ and $\u A$ be as in \S\ref{Se:Gmu-groupoid}.
We will need a Lie version and an ideal version of the $(G,\mu)$-groupoids of \S\ref{Se:displays-Gmu-groupoid}.
This is mainly a collection of notation.
Let $\Ad:G\to\uAut(G)$ and $\Ad:G\to\GL(\Fg)$ denote the left conjugation and   the resulting adjoint action.

\subsection{Ideal version}
\label{Se:Variations-ideal-version}

Let $K$ be an ideal of $\u A$.
Let us write
\beqn
\label{Eq:GK0}
G(K_0)=\ker(G(A_0)\to G(A_0/K_0)),
\eeqn
\beqn
\label{Eq:GKmu}
G(K)_\mu=\ker(G(A)_\mu\to G(A/K)_\mu).
\eeqn
The ring homomorphisms $\sigma,\tau:A\to A_0$ induce group homomorphisms
\beqn
\label{Eq:sigma-tau-GKmu-GK0}
\sigma,\tau:G(K)_\mu\to G(K_0).
\eeqn
For an element $g\in G(A_0)$ we define a group homomorphism 
\beqn
\sigma_g:G(K)_\mu\to G(K_0),\quad \sigma_g=\Ad(g)\circ\sigma.
\eeqn
The twisted conjugation groupoid for $(\sigma_g,\tau)$ will be denoted by
\beqn
\label{Eq:Gamma-GmuKg}
\pWin^{G,\mu}(\u A,K,g)=(G(K_0)\tcg{\sigma_g}{\tau}G(K)_\mu).
\eeqn
The associated orbit map of the unit element is the map
\beqn
\label{Eq:gamma-g}
\gamma_g:G(K_\mu)\to G(K_0),\quad\gamma_g(x)=\tau(x)^{-1}\sigma_g(x).
\eeqn

\subsection{Equivariance}

The assignment $g\mapsto\pWin^{G,\mu}(\u A,K,g)$ extends to a functor
\beqn
\label{Eq:GammaGmuAK-functor}
\pWin^{G,\mu}(\u A,K,-):\pWin^{G,\mu}(\u A)\to(\mathrm{groupoids})
\eeqn
as follows. An element $h\in G(K)_\mu$ gives a morphism $g*h\to g$ in $\pWin^{G,\mu}(\u A)$, to which we associate the isomorphism of groupoids
\beqn
\pWin^{G,\mu}(\u A,K,g*h)\xrightarrow\sim\pWin^{G,\mu}(\u A,K,g)
\eeqn
given by $\Ad(h)$ on $G(K)_\mu$ and by $\Ad(\tau(h))$ on $G(K_0)$. 

\subsection{Lie version}
\label{Se:variations-Lie-version}

The given action of $\GG_m$ on $G_{W(k)}$ (by right conjugation with $\mu$) 
induces an action of $\GG_m$ on the vector group $\Fg_{W(k)}$. Thus we have $A_0$-modules $\Fg(A_0)$ and $\Fg(A)_\mu$, explicitly
$\Fg(A_0)=\Fg\otimes_{\ZZ_p}A_0$ and 
$\Fg(A)_\mu=\bigoplus_{i}\Fg^i\otimes_{W(k)}A_i$
where $\Fg^i$ is given by \eqref{Eq:grading-Fg-Wk}.
The ring homomorphisms $\sigma,\tau:A\to A_0$ give group homomorphisms
$\sigma,\tau:\Fg(A)_\mu\to\Fg(A_0)$,
and for $g\in G(A_0)$ we consider the twisted conjugation groupoid
\beqn
\label{Eq:GammagmuAg}
\pWin^{\Fg,\mu}(\u A,g)=(\Fg(A_0)\tcg{\sigma_g}{\tau} \Fg(A)_\mu).
\eeqn
where $\sigma_g=\Ad(g)\circ\sigma$ as before (with a different meaning of $\Ad$).
Since $\Fg$ is commutative, $\pWin^{\Fg,\mu}(\u A,g)$ is the underlying groupoid of the Picard category associated to the complex of abelian groups placed in degrees $\{-1,0\}$
\beqn
\label{Eq:Cgmu-uAg}
C^{\Fg,\mu}(\u A,g)=[\Fg(A)_\mu\xrightarrow{\sigma_g-\tau}\Fg(A_0)].
\eeqn

\subsection{Adjoint version}
\label{Se:variations.adjoint}

Let us reformulate \eqref{Eq:Cgmu-uAg}.
The composition $\tilde\mu={\Ad}\circ\mu$ is a cocharacter of $\GL(\Fg)$ defined over $W(k)$.
The homomorphism $\Ad:G\to\GL(\Fg)$ induces a functor of groupoids
$\Ad^{G,\mu}:\pWin^{G,\mu}(\u A)\to
\pWin^{\GL(\Fg),\tilde\mu}(\u A)$
and thus a functor
\beqn
\label{Eq:tildeAdGmu}
\Lie^{G,\mu}:\pWin^{G,\mu}(\u A)\to\Win(\u A),\quad 
\Lie^{G,\mu}=\real\circ\Ad^{G,\mu}
\eeqn
using the functor $\real$ of Lemma \ref{Le:GLn-groupoid}. Explicitly, $\Lie^{G,\mu}(g)=(\Fg_{\tilde\mu}\otimes_{W(k)}A,\Ad(g))$ where $\Fg_{\tilde\mu}$ is the graded $W(k)$-module of \S\ref{Se:Groupoids.GL(V)} for $V=\Fg$.
It follows that for $g\in G(A_0)$
\beqn
\label{Eq:CgmuAg-C0}
C^{\Fg,\mu}(\u A,g)\cong \Gamma_0(\Lie^{G,\mu}(g))
\eeqn
using the complex $\Gamma_0(-)$ of \eqref{Eq:Ci}.
Indeed, we have $\Fg_{\tilde\mu}=\bigoplus_{i}\Fg^i(i)$
because $\Fg_{\tilde\mu,i}=\Fg^{-i}$ (the inversion of degrees appears since $\Fg^i$ refers to the right conjugation action of $\mu$, while $\Ad$ is defined by left conjugation) and thus $(\Fg_{\tilde\mu,A})_0=\Fg(A)_\mu$.

\subsection{Adding coefficients}

\label{Se:variations.adjoint-coeff}

For an $\u A$-module $\u N=(N,F_N)$ and $g\in G(A_0)$ we define a complex
\beqn
\label{Eq:CgmuNg}
C^{\Fg,\mu}(\u N,g)=\Gamma_0(\Lie^{G,\mu}(g)\otimes\u N)
\eeqn
and denote by $\pWin^{\Fg,\mu}(\u N,g)$ the associated quotient groupoid. 
The assignments $g\mapsto C^{\Fg,\mu}(\u N,g)$ and $g\mapsto \pWin^{\Fg,\mu}(\u N,g)$ for varying $g$ extend to functors
\beqn
\label{Eq:CgmuN-functor}
C^{\Fg,\mu}(\u N,-):\pWin^{G,\mu}(\u A)\to Ch(\Ab),
\eeqn
\beqn
\label{Eq:GammagmuN-functor}
\pWin^{\Fg,\mu}(\u N,-):\pWin^{G,\mu}(\u A)\to (\mathrm{groupoids})
\eeqn
since $\Lie^{G,\mu}$ is a functor.
The notation does not record the frame $\u A$, which is justified by:

\bRe
\label{Re:CgmuNg-CgmuN'g'}
If $\u A\to \u A'$ is a homomorphism of frames and $\u N'$ is an $\u A'$-module, hence also an $\u A$-module which we denote by $\u N$, then for $g\in G(A_0)$ with image $g'\in G(A_0')$ we have a natural isomorphism of complexes
$\label{Eq:CgmuNg-cGmuNg'}
C^{\Fg,\mu}(\u N,g)\cong C^{\Fg,\mu}(\u N',g')$.
\eRe

\subsection{Exponential}

An ideal $K$ of $\u A$ gives an $\u A$-module $\u K$ as in Remark \ref{Re:ideal-module}, and we have
\beqn
C^{\Fg,\mu}(\u K,g)=[\Fg(K)_\mu\xrightarrow{\sigma_g-\tau}\Fg(K_0)]
\eeqn
for $g\in G(A_0)$ as a subcomplex of \eqref{Eq:Cgmu-uAg} with 
$\Fg(K_0)=\Fg\otimes_{\ZZ_p}K_0$ and 
$\Fg(K)_\mu=\bigoplus_{i\in\ZZ}\Fg^i\otimes_{W(k)}K_i$.
If $K$ is an ideal of square zero, there are natural isomorphisms of groups
\beqn
\exp:\Fg(K_0)\xrightarrow\sim G(K_0),\qquad
\exp_\mu:\Fg(K)_\mu\xrightarrow\sim G(K)_\mu
\eeqn
which give an isomorphism of groupoids
\beqn
\label{Eq:exp-Gammagmu-GammaGmu}
\exp:\pWin^{\Fg,\mu}(\u K,g)\xrightarrow\sim
\pWin^{G,\mu}(\u A,K,g)
\eeqn
for each $g\in G(K_0)$,
and for varying $g$ an isomorphism of functors on $\pWin^{G,\mu}(\u A)$
\beqn
\label{Eq:exp-groupoid}
\exp:\pWin^{\Fg,\mu}(\u K,-)\xrightarrow\sim\pWin^{G,\mu}(\u A,K,-).
\eeqn

\subsection{Passing to quotients}

Let $J$ be an ideal of the frame $\u A$ and let $\pi:\pWin^{G,\mu}(\u A)\to\pWin^{G,\mu}(\u A/J)$ be the projection. If the $\u A$-module $\u N$ is annihilated by $J$, then $\u N$ is an $\u A/J$-module und thus gives a functor
\beqn
\label{Eq:CgmuN-functor-quot}
C^{\Fg,\mu}(\u N,-):\pWin^{G,\mu}(\u A/J)\to Ch(\Ab)
\eeqn
whose composition with $\pi$ gives back \eqref{Eq:CgmuN-functor}; see Remark \ref{Re:CgmuNg-CgmuN'g'}. Similarly, if $K$ is an ideal of $\u A$ with $JK=0$, then we have a functor
\beqn
\label{Eq:GammaGmuAK-functor-quot}
\pWin^{G,\mu}(\u A,K,-)':\pWin^{G,\mu}(\u A/J)\to(\mathrm{groupoids})
\eeqn
whose composition with $\pi$ gives back \eqref{Eq:GammaGmuAK-functor},
using Lemma \ref{Le:Ad'} below. If the ideal $K$ has square zero, the isomorphism $\exp$ of \eqref{Eq:exp-groupoid} extends to the quotient version of the functors.

\bLe
\label{Le:Ad'}
Let $S$ be a\/ $\ZZ_p$-algebra and let $I,J\subseteq S$ be ideals with $IJ=0$. There is a natural group homomorphism
$\Ad':G(S/J)\to\Aut(G(I))$
such that for $\pi:G(S)\to G(S/J)$ we have ${\Ad'}\circ\pi=\Ad$.
\eLe

\bproof
We use the ring $S\oplus I$ with multiplication $(x,a)(y,b)=(xy,xb+ay+ab)$. The sequence of ring homomorphism $S\leftarrow S\oplus I\to (S/J)\oplus I$ with compatible ideals $I$ in all rings gives the same group $G(I)$ in all three cases, 
and we can define $\Ad'$ by
\[
G(S/J)\to G((S/J)\oplus I)\xrightarrow{\;\Ad\;}\Aut(G(I)).
\qedhere
\]
\eproof

\section{Deformations in $(G,\mu)$-groupoids}
\label{Se:Def-Gmu-groupoids}

Let us recall the unique lifting lemma for $(G,\mu)$-displays and its proof following \cite[\S7.1]{Lau:Higher}.
Let $(G,\mu)$ and $\u A$ be as in \S\ref{Se:Gmu-groupoid} where $\mu$ is $1$-bounded (\S \ref{Se:displays-decompositions}). Let $K$ be an ideal of $\u A$ such that $K_0$ is $p$-complete and $K_0\subseteq tA_1$. This gives a frame homomorphism
\beqn
\label{Eq:pi-uA-uA'}
\pi:\u A\to\u A'=\u A/K
\eeqn
such that the induced ring homomorphism $A_0/tA_1\to A_0'/tA_1'$ is bijective.

\bLe
\label{Le:Gpi0-Gpimu-surjective}
In this situation, the frame homomorphism $\pi$ induces surjective 
group homomorphisms 
$G(\pi_0):G(A_0)\to G(A_0')$ and
$G(\pi)_\mu:G(A)_\mu\to G(A')_\mu$.
\eLe

\bproof
See \cite[Lemma 7.1.4 (b)]{Lau:Higher} and its proof, where the ideal $K$ is assumed to be leveled, but this is not essential. The relation $\sigma_{-1}(t)=p$ implies that $\sigma(tA_0)\subseteq pA_0$ and hence $x^p\in pA_0$ for $x\in tA_0$, thus $x^{p+1}\in pK_0$ for $x\in K_0$, so $K_0/pK_0$ is a nil-ideal. Now one proceeds as in loc.cit.
\eproof

\bLe
\label{Le:functor-pi-G-mu}
The functor of $(G,\mu)$-groupoids
$\pi^{G,\mu}:\pWin^{G,\mu}(\u A)\to \pWin^{G,\mu}(\u A')$
induced by $\pi$ is full (an equivalence) iff for each $g\in G(A_0)$ the map 
\beqn
\label{Eq:gamma-g-spec}
\gamma_g:G(K)_\mu\to G(K_0),\quad \gamma_g(x)=\tau(x)^{-1}\sigma_g(x)
\eeqn
is surjective (bijective). Here $\sigma_g=\Ad(g)\circ \sigma$ as earlier.
\eLe

\bproof
We use Lemma \ref{Le:Gpi0-Gpimu-surjective}. The functor is surjective on objects since $G(\pi_0)$ is surjective, and the functor is fibered in groupoids since $G(\pi)_\mu$ is surjective. Hence the properties `full' and `fully faithful' can be verified on the $1$-categorical fibers. For $g\in G(A)$ with image $g'\in G(A')$, the fiber of $\pi^{G,\mu}$ over $g'$ is isomorphic to the groupoid $\pWin^{G,\mu}(\u A,K,g)$ of \eqref{Eq:Gamma-GmuKg}. This groupoid is connected (equivalent to $*$) iff the orbit map $\gamma_g$ is surjective (bijective).
\eproof

\bRe
Assume that the ideal $K$ is leveled, i.e.\ $\tau_1:K_1\to K_0$ is bijective. Then the group homomorphism $
\tau:G(K)_\mu\to G(K_0)$ is bijective by \cite[Lemma 7.1.4 (a)]{Lau:Higher} (using that $\mu$ is $1$-bounded), and we can define a group homomorphism
\beqn
\dot\sigma=\sigma\circ\tau^{-1}:G(K_0)\to G(K_0).
\eeqn
For $g\in G(K_0)$ let $\dot\sigma_g=\Ad(g)\circ\dot\sigma:G(K_0)\to G(K_0)$.
Then $\gamma_g=\dot\gamma_g\circ\tau$ with
\beqn
\label{Eq:dot-gamma-g}
\dot\gamma_g:G(K_0)\to G(K_0),\quad \dot\gamma_g(x)=x^{-1}\dot\sigma_g(x),
\eeqn
and in Lemma \ref{Le:functor-pi-G-mu} we can replace $\gamma_g$ by $\dot\gamma_g$.
If the endomorphism $\dot\sigma_g$ of $G(K_0)$ is locally nilpotent, $\dot\gamma_g$ is bijective with the following inverse:
For $y\in G(K_0)$ with $\dot\sigma_g^m(y)=0$ we have
\beqn
\label{Eq:dot-gamma-g-inverse}
\dot\gamma_g^{-1}(y)=\dot\sigma_g^{m-1}(y^{-1})\cdots\dot\sigma_g(y^{-1})\cdot y^{-1}.
\eeqn
The following lemma gives a criterion when this holds.
\eRe

\bLe
\label{Le:dot-sigma-loc-nilpotent}
Assume that the ideal $K$ is leveled and $\dot\sigma$-nilpotent (Definition \ref{De:leveled}),
and annihilated by a power of $p$. Then for each $g\in G(A_0)$ the endomorphism $\dot\sigma_g=\Ad(g)\circ\dot\sigma=\sigma_g\circ\tau^{-1}$
of the group $G(K_0)$ is locally nilpotent.
\eLe

\bproof
This is proved in \cite[Lemma 7.1.6]{Lau:Higher}.
\eproof

\subsubsection*{Unique lifting lemma}
Under the assumptions of Lemma \ref{Le:dot-sigma-loc-nilpotent} it follows that for each $g\in G(K_0)$ the map $\dot\gamma_g:G(K_0)\to G(K_0)$ is bijective, thus $\gamma_g$ is bijective, and the functor $\pi^{G,\mu}$ is an equivalence by Lemma \ref{Le:functor-pi-G-mu}. This is a presheaf version of the unique lifting lemma \cite[Proposition 7.1.5]{Lau:Higher}.


\section{Witt and crystalline frames}
\label{Se:Witt-crys-frames}

\subsection{A grading formalism}
\label{Se:frames.grading}

In the following, we will often denote a $\ZZ$-graded abelian group by $M^\oplus=\bigoplus_iM^\oplus_i$ with $M=M^\oplus_0$. 
This will be used in particular if $M$ appears earlier than  $M^\oplus$. The notation $N=\bigoplus_iN_i$ will also be used, but the two notations will not be mixed.

\subsection{Witt frames}

For each $p$-complete ring $R$ there is the Witt frame
\[
\u W(R)=(W(R)^\oplus,\sigma,\tau)
\]
defined in \cite[Example 2.1.3]{Lau:Higher}, using here the notation of \cite{Daniels:A-Tannakian}. This frame is determined by the following relations. We have $W(R)^\oplus_0=W(R)$ as a ring, and $\sigma_0=F:W(R)\to W(R)$ is the Witt vector Frobenius; for $i\ge 1$ we have
\[
W(R)^\oplus_i=I(R)=\ker(W(R)\to R)
\]
as a $W(R)$-module, $\sigma_i:I(R)\to W(R)$ is the inverse of the Verschiebung $V$, the multiplication $t:W(R)^\oplus_i\to W(R)^\oplus_{i-1}$ is the inclusion $I(R)\to W(R)$ when $i=1$ and $p:I(R)\to I(R)$ when $i\ge 2$;
the multiplication of homogeneous elements of positive degrees in $W(R)^\oplus$ is given by $V(a)\cdot V(b)=V(ab)$ for $a,b\in W(R)$.
See \cite[10.2.6]{Drinfeld:The-Lau} for an alternative description.

\subsubsection*{Truncated version}

If $R$ is an $\FF_p$-algebra, we have the $n$-truncated Witt frame
\beqn
\label{Eq:uWnR}
\u{W\!}_n(R)=(W_n(R)^\oplus,\sigma,\tau)
\eeqn
as a quotient of $\u W(R)$ as defined in \cite[Example 2.1.6]{Lau:Higher}. Here $W_n(R)^\oplus_0=W_n(R)$ is the ring of $n$-truncated Witt vectors, and $W_n(R)^\oplus_i=I_{n+1}=\ker(W_{n+1}(R)\to R)$ for $i\ge 1$.

\bRe
If $R$ is a semiperfect $\FF_p$-algebra, which means that the Frobenius map $R\to R$ is surjective, then $\u{W\!}_n(R)=\u W(R)/p^n$.
\eRe

\bRe
As a functor of the $\FF_p$-algebra $R$, the components of the graded rings $W(R)^\oplus$ and $W_n(R)^\oplus$ are representable by affine spaces (of infinite or finite dimension). Thus $\u W(-)$ and $\u {W\!}_n(-)$ define frame schemes and in particular fpqc sheaves of frames on $\Aff_{\FF_p}$.
\eRe

\subsection{Displays and truncated displays} 

A display over a $p$-complete $\ZZ_p$-algebra $R$ is a window over the frame $\u W(R)$, and an $n$-truncated display over an $\FF_p$-algebra $R$ is a window over the frame $\u{W\!}_n(R)$. These objects form additive categories 
\[
\Disp(R)=\Win(\u W(R)), \qquad
\Disp_n(R)=\Win(\u {W\!}_n(R)).
\]
We denote by $\Disp^{\eff}(R)$ and $\Disp^{\eff}_n(R)$ the full subcategories of effective (truncated) displays. 

\subsection{Some ideals in Witt frames}
\label{Se:Witt-crys-frames.ideals-Witt}

If $R\to R'$ is a homomorphism of $\ZZ_p$-algebras with kernel $\Fb$, let $W(\Fb)^\oplus$ denote the kernel of the frame homomorphism $\u W(R)\to\u W(R')$. If $p$ is nilpotent in $R$, there is an ideal $\hat W(R)\subseteq W(R)$ which consists of all Witt vectors $x=(x_0,x_1,\ldots)$ such that all $x_i$ are nilpotent and only finitely many $x_i$ are non-zero. This extends to an ideal $\hat W(R)^\oplus$ of the frame $\u W(R)$ such that for $m>0$, the isomorphisms
\[
\tau_{-m}: W(R)^\oplus_{-m}\xrightarrow\sim W(R),\qquad
\sigma_m: W(R)^\oplus_{m}\xrightarrow\sim  W(R)
\]
restrict to isomorphisms
\[
\tau_{-m}:\hat W(R)^\oplus_{-m}\xrightarrow\sim\hat W(R),\qquad
\sigma_m:\hat W(R)^\oplus_{m}\to\hat W(R).
\]
Moreover let $\hat W(\Fb)^\oplus=W(\Fb)^\oplus\cap \hat W(R)^\oplus$; this is an ideal of the frame $\u W(R)$. We will denote by $\u{\hat W}(\Fb)$ the associated $\u W(R)$-module as in Remark \ref{Re:ideal-module}.

\bRe
\label{Re:Wb-module-WnR}
If $\Fb$ is annihilated by the $n$th Frobenius power $F^n$, then the ideal $W(\Fb)^\oplus$ is annihilated by the kernel of $W(R)^\oplus\to W_n(R)^\oplus$, and hence $\u{\hat W}(\Fb)$ is a module for $\u {W\!}_n(R)$. This follows from the relation $V^n(a)b=V^n(aF^n(b))$ for $a,b\in W(R)$.
\eRe

\subsection{Crystalline frames}

Let $R$ be a quasiregular semiperfect $\FF_p$-algebra; see \cite[Definition 8.8]{BMS:Topological}. Then $A_{\crys}(R)$ is torsion-free by \cite[Theorem 8.14 (1)]{BMS:Topological}, and there is a frame
\[
\u A_{\crys}(R)=(A_{\crys}(R)^\oplus,\sigma,\tau)
\]
where $A_{\crys}(R)^\oplus$ is the Rees ring of the Nygaard filtration. Explicitly, if $\sigma_0=\varphi$ denotes the Frobenius of $A_{\crys}(R)$, the degree $i$ component is given by
\[
A_{\crys}(R)^\oplus_i=\NNN^{\ge i}A_{\crys}(R)=\{x\in A_{\crys}(R)\mid\varphi(x)\in p^iA_{\crys}(R)\}
\]
with $\sigma_i=p^{-i}\varphi$ on $\NNN^{\ge i}A_{\crys}(R)$. The evident inclusion $I_{\crys}(R)\subseteq\NNN^{\ge 1}A_{\crys}(R)$ is an equality by \cite[Theorem 8.14 (4), Proposition 8.12]{BMS:Topological}. The frame $\u A_{\crys}(R)$ is $p$-complete. The reduction modulo $p^n$ of the frame $\u A_{\crys}(R)$ will be denoted by
\[
\u A_n(R)=\u A_{\crys}(R)/p^n=(A_n(R)^\oplus,\sigma,\tau).
\]
The category $\Aff_{\FF_p}^{\qrsp}$ of quasiregular semiperfect affine $\FF_p$-schemes becomes a site using the quasi-syntomic topology.
We have presheaves of frames $\u A_{\crys}$ and $\u A_n$ on this category.

\bLe
\label{Le:frame-Acrys-qsyn-sheaf}
The presheaves $\u A_{\crys}$ and $\u A_n$ are quasi-syntomic sheaves on $\Aff_{\FF_p}^{\qrsp}$.
\eLe

\bproof
This can be extracted from the results of \cite[\S 8]{BMS:Topological}. Namely, by \cite[Proposition 8.12]{BMS:Topological} there is an isomorphism
$A_1(R)=A_{\crys}(R)/p\cong L\Omega_{R/\FF_p}$ (derived de Rham cohomology), and $A_1(R)=\colim_i{\Fil_i^{\conj}}A_1(R)$ with 
${\Fil_i^{\conj}}A_1(R)\cong{\Fil_i^{\conj}}L\Omega_{R/\FF_p}$
(conjugate filtration), moreover ${\gr_i^{\conj}}A_1(R)\cong\bigwedge^iL_{R/\FF_p}[-i]$ by the derived Cartier isomorphism, and $\bigwedge^nL_{R/\FF_p}$ is an fpqc sheaf by \cite[Proposition 3.1]{BMS:Topological}. Hence ${\Fil_i^{\conj}}A_1$ and then also $A_1$ and all $A_n$ and $A_{\crys}$ are quasi-syntomic sheaves. By \cite[Proposition 8.14]{BMS:Topological}, for $i\ge 0$ there is an exact sequence $0\to\NNN^{\ge{i+1}}A_{\crys}(R)\to\NNN^{\ge i}A_{\crys}(R)\to F_i\to 0$ with $F_i={\Fil_i^{\conj}}A_1$, which gives an exact sequence of presheaves 
\[
0\to F_i\to(\NNN^{\ge{i+1}}A_{\crys})/p\to(\NNN^{\ge{i}}A_{\crys})/p\to F_i\to 0.
\] 
By induction it follows that $(\NNN^{\ge{i}}A_{\crys})/p$ is a sheaf; then $\NNN^{\ge{i}}A_{\crys}$ is a sheaf as well.
\eproof

\subsection{Comparison map}

Let $R$ by a quasiregular semiperfect $\FF_p$-algebra.

\bLe
\label{Le:rho-exists}
There is a natural homomorphism of frames 
$\varrho:\u A_{\crys}(R)\to\u W(R)$
which induces the identity of $R$. The ring homomorphism $\varrho$ is surjective.
\eLe


\bproof 
By the universal property of $W(R)$ there is a unique homomorphism of $R$-augmented $\delta$-rings $\rho:A_{\crys}(R)\to W(R)$. We claim that $\rho$ is a pd homomorphism and a homomorphism of $1$-frames, i.e.\ $\rho$ commutes with the divided Frobenius. Indeed, $\rho$ is the composition of homomorphisms of $\delta$-rings
$A_{\crys}(R)\to W(A_{\crys}(R))\to W(R)$
where the first arrow preserves the divided powers and the divided Frobenius since $W(A_{\crys}(R))$ is torsion-free, and the second arrow preserves the divided powers and the divided Frobenius by the functoriality of Witt vectors. Now $\rho$ extends to a unique frame homomorphism $\varrho$ with $\varrho_0=\rho$ because for $i>0$ the component $\varrho_i$ is determined by the following commutative diagram.
\[
\xymatrix@M+0.2em{
\NNN^{\ge i}A_{\crys}(R) \ar[r]^-{\varrho_i} \ar[d]_{\sigma_i} &
W(R)^\oplus_i \ar[d]_{\sigma_i} \ar@{=}[r] & I(R) \ar[dl]^{V^{-1}}_\sim \\
A_{\crys}(R) \ar[r]^-{\varrho_0} &
W(R)
}
\] 
In other words $\varrho_i(x)=V(\varrho_0(\sigma_i(x)))$ as an element of $I(R)$. The verification that $\varrho$ is a frame homomorphism is straightforward (using that $\rho$ preserves $\dot\sigma$). Clearly $\varrho_0$ is surjective, then $\varrho_1$ is surjective, and $\varrho_i$ is surjective for $i\ge 2$ since $\varrho_1$ factors as
\[
\NNN^{\ge 1}A_{\crys}(R)\xrightarrow{p^{i-1}}\NNN^{\ge 1}A_{\crys}(R)\xrightarrow{\varrho_i}I(R).
\qedhere
\] 
\eproof

\subsubsection*{Truncated version}
The homomorphism $\varrho$ induces a surjective frame homomorphism
\beqn
\label{Eq:rhon}
\varrho_n:\u A_n(R)\to\u{W\!}_n(R).
\eeqn
Here $\varrho_n$ is the reduction of $\varrho$, not the degree $n$ component of $\varrho$ as in the proof of Lemma \ref{Le:rho-exists}.

\section{Witt and crystalline $(G,\mu)$-displays}
\label{Se:Witt-crys-Gmu-disp}

Let $(G,\mu)$ be a pair as in \S\ref{Se:Gmu-groupoid}. Later we assume that $\mu$ is $1$-bounded.

\subsection{Witt vector $(G,\mu)$-displays}
\label{Se:Witt-crys.Witt}

If $R$ is a $k$-algebra, the frame $\u {W\!}_n(R)$ of \eqref{Eq:uWnR} is a frame over $W(k)$, and we can form the $(G,\mu)$-groupoid
\beqn
\label{GammaGmuWnR}
\pWin^{G,\mu}(\u{W\!}_n(R))=(G(W_n(R))\tcg{\sigma}{\tau}G(W_n(R)^\oplus)_\mu)
\eeqn
as in \eqref{Eq:Gamma-GmuA}. For varying $R$ we obtain a presheaf of groupoids $\pWin^{G,\mu}(\u {W\!}_n(-))$ over $\Aff_k$, which is represented by a smooth groupoid of smooth affine $k$-schemes. The associated algebraic stack will be denoted by
\[
\Disp^{G,\mu}_n
=[\pWin^{G,\mu}(\u {W\!}_n(-))]_{\et};
\]
this is the stack of $(G,\mu)$-displays of level $n$, which is a smooth algebraic stack over $k$ of pure dimension zero with affine diagonal; see \cite[Example 5.3.7]{Lau:Higher}. It follows that $\Disp_n^{G,\mu}$ is an fpqc stack and hence equal to 
$[\pWin^{G,\mu}(\u {W\!}_n(-))]_\tau$
for any topology $\tau$ between etale and fpqc.

\subsection{Crystalline $(G,\mu)$-displays}
\label{Se:Witt-crys.crys}

If $R$ is a quasiregular semiperfect $k$-algebra, the frame $\u A_n(R)=\u A_{\crys}(R)/p^n$ is a frame over $W(k)$, and we can form the $(G,\mu)$-groupoid
\beqn
\label{GammaGmuAnR}
\pWin^{G,\mu}(\u A_n(R))=(G(A_n(R))\tcg{\sigma}{\tau}G(A_n(R)^\oplus)_\mu).
\eeqn
For varying $R$ we obtain a presheaf of groupoids $\pWin^{G,\mu}(\u A_n(-))$ over $\Aff_k^{\qrsp}$. The associated etale stack will be called the stack of crystalline $(G,\mu)$-displays over $\Aff^{\qrsp}_k$; see \S\ref{Se:Geom-appl} for more.

\bLe
\label{Le:Gammagmu-An-qsyn-sheaf}
The presheaf of groupoids $\pWin^{G,\mu}(\u A_n(-))$ is a quasi-syntomic sheaf of groupoids in the $1$-categorical sense, i.e.\ the presheaves of groups $G(A_n(-))$ and $G(A_n(-)^\oplus)_\mu$ are sheaves.
\eLe

\bproof
By \eqref{Eq:GAmu-decomposition}, \eqref{Eq:Pamu-PA0}, and \eqref{UAmu-uAi} it suffices to show that the presheaves $G(A_n(-))$, $P^-(A_n(-))$, and $\Fg^i\otimes_{W(k)} A_n(-)^\oplus_i$ for $i>0$ are sheaves. This holds since $A_n(-)^\oplus$ is a quasi-syntomic sheaf of graded rings by Lemma \ref{Le:frame-Acrys-qsyn-sheaf}.
\eproof

\subsection{Comparison functor}

For a quasiregular semiperfect $k$-algebra $R$, the frame homomorphism $\varrho_n:\u A_n(R)\to\u {W\!}_n(R)$ of \eqref{Eq:rhon} induces a homomorphism of groupoids
\beqn
\label{Eq:rhon-Gamma}
\varrho_{n}^{G,\mu}:\pWin^{G,\mu}(\u{A}_n(R))\to\pWin^{G,\mu}(\u{W\!}_n(R)).
\eeqn
For $R\in\Aff^{\qrsp}_k$ and $g\in G(A_n(R))$ we denote by $E^{G,\mu}(R,g)$ be the inertia of $\varrho_{n}^{G,\mu}$ at $g$, thus
\beqn
\label{Eq:EGmunRg}
E^{G,\mu}_n(R,g)=\{x\in G(A_n(R))_\mu\mid g*x=g, \;\varrho_n(x)=1\text{ in }G(W_n(R))\}.
\eeqn
This group is functorial in $g$ by conjugation: Each $h\in G(A_n(R)^\oplus)_\mu$ gives a morphism $g*h\to g$ in $\pWin^{G,\mu}(\u A_n(R))$, which acts by
$\Ad(h):E_n^{G,\mu}(g*h)\xrightarrow\sim E_n^{G,\mu}(g)$.
This gives a functor
\beqn
E^{G,\mu}(R,-):\pWin^{G,\mu}(\u A_n(R))\to(\mathrm{groups}).
\eeqn
The following is the central technical result of this article, which will be proved in \S \ref{Se:Proof-Prop}.

\bPr
\label{Pr:Main}
Assume that $\mu$ is $1$-bounded.
The functor $\varrho_n^{G,\mu}$ of \eqref{Eq:rhon-Gamma} is surjective on objects and full, and for $R\in\Aff^{\qrsp}_k$ and $g\in G(A_n(R))$ with image $\bar g\in G(W_n(R))$ there is an isomorphism of groups
\beqn
\label{Eq:Pr-Main}
E_n^{G,\mu}(R,g)\cong H^0(C^{\Fg,\mu}(\u{\hat W}(R[F^n]),\bar g))
\eeqn
which is functorial in $R$ and in $g$. In particular, the group $E_n^{G,\mu}(R,g)$ is abelian.
\ePr

Here $C^{\Fg,\mu}(-,-)$ is the complex of abelian groups in degree $\{-1,0\}$ of \eqref{Eq:CgmuNg}, explicitly
\beqn
\label{Eq:ZngmuRg}
C^{\Fg,\mu}(\u{\hat W}(R[F^n]),g)=
\Gamma_0(\Lie^{G,\mu}(g)\otimes\u{\hat W}(R[F^n])).
\eeqn
We will also see that $H^{-1}(C^{\Fg,\mu}(\u{\hat W}(R[F^n]),\bar g))$ is trivial.

\section{Sheared Witt vectors (semiperfect case)}
\label{Se:Sheared}

As a technical device, for a semiperfect ring $R$ we will use a  modification $\sW(R)$ of the ring of Witt vectors $W(R)$. A generalisation that applies to every $p$-nilpotent ring $R$ is studied in \cite{{Dr25:Ring-stacks-conjecturally},BMVZ}. For consistency we adopt the notation $\sW(R)$ and the terminology \emph{sheared Witt vectors} of loc.cit.; an earlier version of this text used the notation $\check W(R)$.

\subsection{Construction and basic properties}

Let $R$ be a semiperfect $\FF_p$-algebra. 
Let $R=R^\flat/J$ where $R^\flat=R^{\perf}$ is the limit perfection. 
Let $J_n=J/F^n(J)$ as an ideal of $R_n=R^\flat/F^n(J)$.
Then $J=\varprojlim_n J_n$. We set
\beqn
\hat W(J)=\varprojlim_n\hat W(J_n),
\eeqn
which is an ideal of $W(R^\flat)$ contained in $W(J)$, and 
\beqn
\sW(R)=W(R^\flat)/\hat W(J).
\eeqn
The ring $\sW(R)$ is derived $p$-complete because $W(R)$ is $p$-complete and each $\hat W(J_n)$ is annihilated by $p^n$, but $\sW$ is classically $p$-complete only when $R$ is perfect since $\sW(R)/p^n=W_n(R)$, so the classical $p$-completion of $\sW(R)$ is $W(R)$.
The ring $\sW(R)$ carries operators $F$ and $V$ induced by those of $W(R^\flat)$. Let 
\beqn
\sI(R)=V(\sW(R))=\ker(\sW(R)\to R)
\eeqn
and let $\dot F=V^{-1}:\sI(R)\to\sW(R)$.
There is an exact sequence
\beqn
\label{Eq:K-sW-W}
0\to K(R)\to\sW(R)\to W(R)\to 0
\eeqn
where $K(R)=W(J)/\hat W(J)$ is an ideal of square zero on which $V$ and $\dot F$ are bijective, $F$ is injective, thus $p$ is injective. The square zero property can be verified directly, or one can invoke \cite[\href{https://stacks.math.columbia.edu/tag/0G3G}{Lemma 0G3G}]{Stacks}.
The exact sequence \eqref{Eq:K-sW-W} is the limit over $n$ of exact sequences
\beqn
0\to K_n(R)\to\sW^{(n)}(R)\to W(R)\to 0
\eeqn
with $K_n(R)=W(J_n)/\hat W(J_n)$ and $\sW^{(n)}(R)=W(R_n)/\hat W(J_n)$.

\subsection{Frame structure}
\label{Se:Sheared-Witt.Frame}

In the definition of the frame $\u W(R)$ we can replace the ring $W(R)$ by $\sW(R)$ or $\sW^{(n)}(R)$. For a short formulation, a frame $\u A$ will be called even if for $i\ge 1$ the homomorphism $\sigma_i:A_i\to A_0$ is bijective. Thus $\u W(R)$ is even. The following is straightforward.

\bLe
\label{Le:frame-sW}
Let $R$ be a semiperfect $\FF_p$-algebra.
There are unique even frames $\u\sW(R)$ and $\u\sW^{(n)}(R)$ with surjective frame homomorphisms $\u W(R^\flat)\to\u\sW(R)\to\u\sW^{(n)}(R)$ which give the projections $W(R^\flat)\to\sW(R)\to\sW^{(n)}(R)$ in degree zero. \qed
\eLe

\subsection{Relation with \boldmath $A_{\crys}(R)$}

Assume now that $R$ is quasiregular semiperfect over $\FF_p$. Let $N(R)$ be the kernel of the natural homomorphism $\rho:A_{\crys}(R)\to W(R)$. We have a $\varphi$-linear endomorphism 
\[
\dot\varphi=p^{-1}\varphi:N(R)\to N(R).
\] 
Let $N_n(R)=N(R)/p^nN(R)$ as an ideal of 
\[
B_n(R)=A_{\crys}(R)/p^nN(R),
\]
let $N_n(R)^{\nil}$ be the set of elements of $N_n(R)$ on which $\dot\varphi$ is nilpotent, and
\[
N(R)^{\nil}=\varprojlim_nN_n(R)^{\nil}.
\]
Let $\varkappa:W(R^\flat)\to A_{\crys}(R)$ be the canonical ring homomorphism. We have $\varkappa\circ F=\varphi\circ\varkappa$.

\bLe
\label{Le:varkappa-n}
The ring homomorphism $\varkappa$ induces $\varkappa^{(n)}:W(R_n)\to B_n(R)$.
\eLe

\bproof
The homomorphism $\varkappa$ restricts to $W(J)\to N(R)$ and then to $W(F^n(J))\to p^nN(R)$ since $\varphi=p\dot\varphi$ on $N(R)$. Passing to the quotients gives $\varkappa^{(n)}$.
\eproof

\bLe
\label{Le:varkappa-restrict}
The ring homomorphisms $\varkappa$ and $\varkappa^{(n)}$ restrict to
\[
\varkappa:\hat W(J)\to N(R)^{\nil},\qquad
\varkappa^{(n)}:\hat W(J_n)\to N_n(R)^{\nil}.
\]
\eLe

\bproof
For $a\in J$, the elements $V^m([a])\in W(J)$ and $\varkappa(V^m[a])\in N(R)$ satisfy
\[
\dot\varphi^m(\varkappa(V^m([a])))=\varkappa([a]),\qquad
\dot\varphi^r(\varkappa([a]))=(p^r-1)!\gamma_{p^r}(\varkappa([a]))
\] 
where $\gamma_{p^r}(\varkappa([a]))\in N(R)$. Indeed, $N(R)$ is a pd ideal of the torsion-free ring $A_{\crys}(R)$, and the relations hold after multiplication by a power of $p$. Since $(p^r-1)!$ is divisible by $p^{r-1}$, it follows that $\varkappa(V^m([a]))$ lies in $N(R)^{\nil}$. This yields the assertion for $\varkappa^{(n)}$, and the assertion for $\varkappa$ follows in the limit.
\eproof

\bLe
\label{Le:ker-Nn-barNn}
Let $\bar N_n(R)$ be the kernel of $A_{\crys}(R)/p^n\to W_n(R)$. There is an exact sequence of $A_{\crys}(R)$-modules
\beqn
\label{Eq:ker-Nn-barNn}
0\longrightarrow W(R[F^n])\xrightarrow{\;\alpha_n\;}N_n(R)\longrightarrow\bar N_n(R)\longrightarrow 0
\eeqn
where $\alpha_n$ is induced by multiplication by $p^n$ in $A_{\crys}(R)$.  \eLe

\bproof
We have a commutative diagram of $A_{\crys}(R)$-modules with exact rows, using that $A_{\crys}(R)$ is torsion-free and $p^n=V^nF^n$ in $W(R)$. 
\beqn
\xymatrix@M+0.2em{
0 \ar[r] & W(R) \ar[r]^-{p^n} \ar[d] & B_n(R) \ar[r] \ar[d] & A_{\crys}(R)/p^n \ar[r] \ar[d] & 0 \\
0 \ar[r] & W(R/R[F^n]) \ar[r]^-{p^n} & W(R) \ar[r] & W_n(R) \ar[r] & 0
}
\eeqn
The vertical arrows are surjective, and the sequence of their kernels gives \eqref{Eq:ker-Nn-barNn}. 
\eproof

\bLe
\label{Le:WJ-NRnil-NR-surj}
The homomorphism $W(J)\oplus N(R)^{\nil}\to N(R)$ is surjective.
\eLe

\bproof
Since all modules in the statement are derived $p$-complete, it suffices to show that $\hat W(J)\oplus N(R)^{\nil}\to N(R)/p$ is surjective; see \cite[\href{https://stacks.math.columbia.edu/tag/09B9}{Lemma 09B9}]{Stacks}. We use the exact sequence \eqref{Eq:ker-Nn-barNn} for $n=1$. Since $A_{\crys}(R)/p$ is an $R^\flat$-algebra generated by the elements $\gamma_m([a])$ for $a\in J$ and $m\ge 1$, which lie in the image of $N(R)^{\nil}$, the homomorphism $N(R)^{\nil}\to\bar N_1(R)$ is surjective. The image of $\alpha_1$ lies in the image of $W(J)\to N_1(R)$. The lemma follows.
\eproof

\bPr
\label{Pr:srho}
There is a unique pd homomorphism $\srho:A_{\crys}(R)\to\sW(R)$ such that the following diagram commutes; the outer arrows are canonical.
\[
\xymatrix@M+0.2em{
W(R^\flat) \ar[r] \ar[d] & \sW(R) \ar[d] \\
A_{\crys}(R) \ar[r]^-\rho \ar@{-->}[ru]^{\srho} & W(R)
}
\] 
The ring homomorphism $\srho$ is a homomorphism of $1$-frames (\S\ref{Se:Frames.Terminology}).
\ePr

\bproof
Let $A=W(R^\flat)$ and let $I$ be the kernel of $A\to R$. Then $A_{\crys}(R)$ is the $p$-completion of the pd envelope $D_I(A)$ relative to the canonical divided powers of $p$.  Here $D_I(A)$ is torsion-free by \cite[Remark 8.3]{Mathew:Some-recent}, and hence $A_{\crys}(R)$ is also the derived $p$-completion of $D_I(A)$. The homomorphism $W(R^\flat)\to\sW(R)$ extends to a unique pd homomorphism $\srho':D_I(A)\to\sW(R)$, which makes the diagram commute with $D_I(A)$ in place of $A_{\crys}(R)$. Taking the derived $p$-completion gives the diagram of the proposition because $\sW(R)$ is derived $p$-complete. The homomorphism $\srho$ is a homomorphism of $1$-frames if $\srho\circ\varphi=F\circ\srho$ and $\srho\circ\dot\varphi=\dot F\circ\srho$ on $I_{\crys}(R)$. It suffices to verify the corresponding relations for $\srho'$. We have $\srho'\circ\varphi=F\circ\srho'$
by the universal property of $D_I(A)$. The ring $D_I(A)$ is generated over $W(R^\flat)$ by the elements $\gamma_m([a])$ for $a\in J$ and $m\ge 1$, and the canonical pd ideal of $D_I(A)$ is generated by these elements and $p$. The relation $\srho'\circ\dot\varphi=\dot F\circ\srho'$ follows since $\srho'$ annihilates all $\gamma_m([a])$, and we have $\dot\varphi(\gamma_m([a]))=((pm!)/(p\cdot m!))\cdot\gamma_{pm}([a])$.
\eproof

\bPr[Theorem \ref{Th:Acris-sW}]
The ring homomorphism $\varkappa$ induces an isomorphism 
\[
\sW(R)\xrightarrow\sim A_{\crys}(R)/N(R)^{\nil}.
\]
\ePr

\bproof
We claim that there is the following commutative diagram with injective vertical maps and Cartesian squares; then the proposition follows.
\beqn
\label{Eq:Diag-sW-Acris}
\xymatrix@M+0.2em{
W(R^\flat) \ar[r]^-{\varkappa} & A_{\crys}(R) \ar[r]^-{\srho} & \sW(R) \ar[r] & W(R) \\
W(J) \ar[r] \ar[u] & N(R) \ar[r] \ar[u] & K(R) \ar[r] \ar[u] & 0 \ar[u] \\
\hat W(J) \ar[r] \ar[u] & N(R)^{\nil} \ar[r] \ar[u] & 0 \ar[u] 
}
\eeqn
Clearly the upper squares exist and are Cartesian, and the kernel of $W(J)\to K(R)$ is $\hat W(J)$. Using Lemmas \ref{Le:varkappa-restrict} and \ref{Le:WJ-NRnil-NR-surj} it suffices to show that the homomorphism $\srho:N(R)\to K(R)$ annihilates $N(R)^{\nil}$. This homomorphism is the limit over $n$ of homomorphisms
\beqn
\label{Eq:srhon-NnR-WW}
\srho_n:N_n(R)\to W(J_n)/\hat W(J_n)
\eeqn
since $W(J_n)$ is annihilated by $F^n$ and hence by $p^n$, moreover we have $\srho_n\circ\dot\varphi=\dot F\circ\srho_n$ by the corresponding relation for $\srho$. Since $\dot F$ is bijective on the target of $\srho_n$, this map annihilates $N_n(R)^{\nil}$, hence $\srho$ annihilates $N(R)^{\nil}$.
\eproof

\bCo
\label{Co:sWnR-BnR/NnRnil}
The ring homomorphism $\varkappa^{(n)}$ of Lemma \ref{Le:varkappa-n} induces an isomorphism 
\[
\sW^{(n)}(R)\cong B_n(R)/N_n(R)^{\nil}.
\]
\eCo

\bproof
We claim that there is the following commutative diagram with injective vertical maps and Cartesian squares; this is a variant of \eqref{Eq:Diag-sW-Acris}.
\[
\xymatrix@M+0.2em{
W(R_n) \ar[r]^-{\varkappa^{(n)}} & B_n(R) \ar[r]^-{\srho^{(n)}} & \sW^{(n)}(R) \ar[r] & W(R) \\
W(J_n) \ar[r] \ar[u] & N_n(R) \ar[r] \ar[u] & K_n(R) \ar[u] \ar[r] & 0 \ar[u] \\
\hat W(J_n) \ar[r] \ar[u] & N_n(R)^{\nil} \ar[r] \ar[u] & 0 \ar[u] 
}
\]
The homomorphism $\srho$ composed with $\sW(R)\to\sW^{(n)}(R)$ restricts to $N(R)\to K_n(R)$ where $K_n(R)$ is annihilated by $F^n$ and hence by $p^n$. Thus $\srho^{(n)}$ exists. Clearly the upper two squares are Cartesian. Now the proof concludes as before.
\eproof

\subsection{Frame structure}

\bPr
\label{Pr:svarrho}
Let $R$ be a quasiregular semiperfect $\FF_p$-algebra.
The ring homomorphism $\srho:A_{\crys}(R)\to\sW(R)$ extends uniquely to a frame homomorphisms $\svarrho:\u A_{\crys}(R)\to\u\sW(R)$, and $\svarrho$ is surjective.
\ePr

\bproof
The homomorphism $\srho$ is a homomorphism of $1$-frames by 
Proposition \ref{Pr:srho}. Since the frame $\u\sW(R)$ is even (see \S\ref{Se:Sheared-Witt.Frame}), the homomorphism $\srho$ extends uniquely to a frame  homomorphism $\svarrho$ where $\svarrho_i:A_{\crys}(R)^\oplus_i\to\sW(R)^\oplus_i$ is given by $\svarrho_i=\sigma_i^{-1}\circ\srho\circ\sigma_i$ when $i\ge 1$. Surjectivity follows as in the proof of Lemma \ref{Le:rho-exists}.
\eproof

\section{An auxiliary frame}
\label{Se:frame-B}

Let $R$ be a quasiregular semiperfect $\FF_p$-algebra.
We will use an auxiliary frame $\u B_n(R)$.

\subsection{The frame}

Let $N(R)^\oplus$ be the kernel of the frame homomorphism $\varrho:\u A_{\crys}(R)\to \u W(R)$ of Lemma \ref{Le:rho-exists}. Then $N(R)^\oplus$ is an ideal of the frame $\u A_{\crys}(R)$.
We define a new frame
\beqn
\u B_n(R)=(B_n(R)^\oplus,\sigma,\tau)=\u A_{\crys}(R)/p^nN(R)^\oplus.
\eeqn
Here $B_n(R)^\oplus$ is a graded ring with $B_n(R)=A_{\crys}(R)/p^nN(R)$ in degree zero.
The frame $\u B_n(R)$ is $p$-complete 
because the ideal $p^nN(R)^\oplus\subseteq A_{\crys}(R)^\oplus$ is $p$-closed in each degree.

\subsection{Relation with \boldmath $\varrho_n$ (plain version)}

There is a commutative diagram of frames for $R$ with surjective homomorphisms:
\beqn
\label{Eq:diagram-frames}
\xymatrix@M+0.2em{
\u A_{\crys}(R) \ar[dr]_\varrho \ar[r]  &
\u B_n(R) \ar[r]^{\pi_n} \ar[d]^{\tilde\varrho_n} &
\u A_n(R) \ar[d]^{\varrho_n}  \\ 
& \u W(R) \ar[r]^-{\bar\pi_n} &
\u {W\!}_n(R)
}
\eeqn
We will denote the kernels of the vertical frame homomorphisms in \eqref{Eq:diagram-frames} by
\beqn
\label{Eq:M-L}
N_n(R)^\oplus=\ker(\tilde\varrho_n)=N(R)^\oplus/p^n N(R)^\oplus,\qquad
\bar N_n(R)^\oplus=\ker(\varrho_n).
\eeqn
As earlier let $R[F^n]$ denote the kernel of $F^n:R\to R$. This gives a graded module $W(R[F^n])^\oplus$ over the graded ring $W_n(R)^\oplus$; see \S\ref{Se:Witt-crys-frames.ideals-Witt} and Remark \ref{Re:Wb-module-WnR}. The following extends Lemma \ref{Le:ker-Nn-barNn}.

\bPr
\label{Pr:ker-tilde-rho-n}
There is an exact sequence of graded $B_n(R)^\oplus$-modules 
\beqn
\label{Eq:ker-tilde-rho-n}
0\longrightarrow W(R[F^n])^\oplus\xrightarrow{\;\alpha_n\;}N_n(R)^\oplus\xrightarrow{\;\pi_n\;}\bar N_n(R)^\oplus\longrightarrow 0
\eeqn
where $\alpha_n$ is induced by multiplication by $p^n$ in $B_n(R)^\oplus$.  Here $W(R[F^n])^\oplus$ is annihilated by the kernel of $B_n(R)^\oplus\to W_n(R)^\oplus$, in particular the image of $\alpha_n$ is an ideal of square zero.
\ePr

\bproof
We have a commutative diagram of $B_n(R)^\oplus$-modules with exact rows, using that $A_{\crys}(R)^\oplus$ is torsion-free and $p^n=V^nF^n$ in $W(R)$.
\beqn
\label{Eq:ker-tilde-rho-n-proof}
\xymatrix@M+0.2em{
0 \ar[r] & W(R)^\oplus \ar[r]^-{p^n} \ar[d]^{\can} & B_n(R)^\oplus \ar[r]^-{\pi_n} \ar[d]^{\tilde\varrho_n} & A_n(R)^\oplus \ar[r] \ar[d]^{\varrho_n} & 0 \\
0 \ar[r] & W(R/R[F^n])^\oplus \ar[r]^-{p^n} & W(R)^\oplus \ar[r]^{\bar\pi_n} & W_n(R)^\oplus \ar[r] & 0
}
\eeqn
The vertical arrows of the diagram are surjective, and the sequence of their kernels gives \eqref{Eq:ker-tilde-rho-n}. The relation $p^n=F^nF^n$ also shows that the $W(R)^\oplus$-module $W(R[F^n])^\oplus$ is annihilated by $p^n$ and is thus a $W_n(R)^\oplus$-module.
\eproof

\bCo
\label{Co:exact-sequence-fiber-prod}
There is an exact sequence
\[
0\longrightarrow
W(R[F^n])^\oplus
\xrightarrow{\;\alpha_n\;}
B_n(R)^\oplus
\xrightarrow{\;\;q\;\;}
A_n(R)^\oplus\times_{W_n(R)^\oplus}W(R)^\oplus
\longrightarrow 0
\]
where $q$ is a homomorphism of graded rings with kernel of square zero.
\eCo

\bproof
The diagram \eqref{Eq:ker-tilde-rho-n-proof} and the exact sequence \eqref{Eq:ker-tilde-rho-n} form a $(3\times 3)$-diagram with exact rows and columns.
\eproof

\subsection{Relation with \boldmath $\varrho_n$ (sheared version)}

In the preceding discussion we replace the frame $\u W(R)$ by its sheared version $\u\sW^{(n)}(R)$ given by Lemma \ref{Le:frame-sW}. To begin with, the frame homomorphism $\svarrho$ of Proposition \ref{Pr:svarrho} induces a frame homomorphism of quotients
\[
{}^s\!\tilde\varrho_n:\u B_n(R)\to\u\sW^{(n)}(R)
\]
because $\u A_{\crys}(R)\xrightarrow{\;\svarrho\;}\u\sW(R)\to\u\sW^{(n)}(R)$ maps $N(R)^\oplus$ into the kernel of $\u\sW^{(n)}(R)\to\u W(R)$, which is annihilated by $F^n$ and hence by $p^n$.
As before we get a commutative diagram of frames for $R$ with surjective homomorphisms:
\beqn
\label{Eq:diagram-frames-s}
\xymatrix@M+0.2em{
\u A_{\crys}(R) \ar[dr]_{\svarrho} \ar[r]  &
\u B_n(R) \ar[r]_{\pi_n} \ar[d]^{{}^s\!\tilde\varrho_n} &
\u A_n(R) \ar[d]^{\varrho_n}  \\ 
& \u \sW^{(n)}(R) \ar[r]_-{{}^s\bar\pi_n} &
\u {W\!}_n(R)
}
\eeqn
Let $\sN_n(R)^\oplus=\ker({}^s\!\tilde\varrho_n)$. 
Then $\sN_n(R)^\oplus_0=N_n(R)^{\nil}$ by Corollary \ref{Co:sWnR-BnR/NnRnil}.

\bPr
\label{Pr:ker-tilde-rho-n-s}
There is an exact sequence of graded $B_n(R)^\oplus$-modules 
\beqn
\label{Eq:ker-tilde-rho-n-s}
0\longrightarrow \hat W(R[F^n])^\oplus\xrightarrow{\;\hat\alpha_n\;}\sN_n(R)^\oplus\xrightarrow{\;\pi_n\;}\bar N_n(R)^\oplus\longrightarrow 0
\eeqn
where $\hat \alpha_n$ is induced by multiplication by $p^n$ in $B_n(R)^\oplus$.  Here $\hat W(R[F^n])^\oplus$ is annihilated by the kernel of $B_n(R)^\oplus\to W_n(R)^\oplus$, in particular the image of $\hat \alpha_n$ is an ideal of square zero.
\ePr

\bproof
In \eqref{Eq:ker-tilde-rho-n-proof} replace $W(R)^\oplus$ by $\sW^{(n)}(R)^\oplus$ and $W(R/R[F^n])$ by $W(R)/\hat W(R[F^n])$ and then proceed as in the proof of Proposition \ref{Pr:ker-tilde-rho-n}. 
\eproof

\bCo
\label{Co:exact-sequence-fiber-prod-s}
There is an exact sequence
\[
0\longrightarrow
\hat W(R[F^n])^\oplus
\xrightarrow{\;\alpha_n\;}
B_n(R)^\oplus
\xrightarrow{\;\;q\;\;}
A_n(R)^\oplus\times_{W_n(R)^\oplus}\sW^{(n)}(R)^\oplus
\longrightarrow 0
\]
where $q$ is a homomorphism of graded rings with kernel of square zero.
\eCo

\bproof
Similar to Corollary \ref{Co:exact-sequence-fiber-prod}.
\eproof

\bLe
The ideals $N_n(R)^\oplus$ and $\sN_n(R)^\oplus$ of the frame $\u B_n(R)$ are leveled, and the ideal $\sN_n(R)^\oplus$ is moreover $\dot\sigma$-nilpotent. See Definition \ref{De:leveled} for the terminology.
\eLe

\bproof
The ideals are leveled because $\u B_n(R)\to\u\sW^{(n)}(R)\to\u W(R)$ are homomorphisms of frames for $R$. By Corollary \ref{Co:sWnR-BnR/NnRnil} the degree zero part of $\sN_n(R)^\oplus$ is $N_n(R)^{\nil}$, which is $\dot\sigma$-nilpotent by definition.
\eproof

\section{Proof of Proposition \ref{Pr:Main}}
\label{Se:Proof-Prop}

Let $(G,\mu)$ be $1$-bounded.
Let $R$ be a quasiregular semiperfect $k$-algebra.
We consider the diagram of frames \eqref{Eq:diagram-frames-s}. 
We have the following group version of Proposition \ref{Pr:ker-tilde-rho-n-s}. 

\bPr
\label{Pr:GhWR-GKnil}
There are exact sequences of groups
\beqn
\label{Eq:GhWR-GKnil}
1  \to
G(\hat W(R[F^n])) \xrightarrow{\;\alpha\;}
G(N_n(R)^{\nil}) \xrightarrow{\;\pi\;}
G(\bar N_n(R)) \to 1
\eeqn
\beqn
\label{Eq:GhWR-GKnil-eins}
1 \to G(\hat W(R[F^n])^\oplus)_\mu \xrightarrow{\;\alpha_\mu\;} 
G(\sN_n(R)^\oplus)_\mu \xrightarrow{\;\pi_\mu\;}
G(\bar N_n(R)^\oplus)_\mu \to 1
\eeqn
where 
$\hat W(R[F^n])^\oplus$ is considered as an ideal of $\u B_n(R)$ by the homomorphism $\hat\alpha_n$ of \eqref{Eq:ker-tilde-rho-n-s}.
\ePr

The graded $W_n(R)^\oplus$-module $\hat W(R[F^n])^\oplus\subseteq W(R[F^n])^\oplus$ is defined in \S\ref{Se:Witt-crys-frames.ideals-Witt}.
The evaluation of $G$ on (homogeneous) ideals is written in \eqref{Eq:GK0} and \eqref{Eq:GKmu}.

\bproof
The functors $G(-_0)$ and $G(-)_\mu$ applied to \eqref{Eq:diagram-frames-s} give commutative diagrams of groups
\[
\xymatrix@M+0.2em{
G(B_n(R)) \ar[r] \ar[d] & G(A_n(R)) \ar[d] & & G(B_n(R)^\oplus)_\mu \ar[r] \ar[d] & G(A_n(R)^\oplus)_\mu \ar[d] \\
G(\sW^{(n)}(R)) \ar[r] & G(W_n(R)) & & G(\sW^{(n)}(R)^\oplus)_\mu \ar[r] & G(W_n(R)^\oplus)_\mu
}
\]
where all homomorphisms are surjective by Lemma \ref{Le:Gpi0-Gpimu-surjective}. The homomorphisms of the kernels of the vertical homomorphisms in these diagrams are the homomorphisms $\pi$ and $\pi_\mu$ of the proposition. 
Since $G$ is smooth and affine, Corollary \ref{Co:exact-sequence-fiber-prod-s} implies that the group homomorphism
\beqn
\label{Eq:GBnR-GAnR-times}
G(B_n(R))\to G(A_n(R))\times_{G(W_n(R))}G(\sW^{(n)}(R))
\eeqn
is surjective with kernel $(\hat W(R[F^n]))$, and, using also \eqref{Eq:GAmu-decomposition} and \eqref{UAmu-uA1}, the group homomorphism
\beqn
\label{Eq:GBnR-GAnR-times-mu}
G(B_n(R)^\oplus)_\mu\to G(A_n(R)^\oplus)_\mu\times_{G(W_n(R)^\oplus)_\mu}G(\sW^{(n)}(R)^\oplus)_\mu
\eeqn
is surjective with kernel $G(\hat W(R[F^n])^\oplus)_\mu$. The proposition follows.
\eproof

\bCo
\label{Co:fundamental-diagram}
For $g\in G(A_n(R))$ with image $\bar g\in G(W_n(R))$ we have the following commutative diagram of pointed sets where all objects are groups, all horizontal arrows are group homomorphisms, the rows are exact, and $\gamma_{g}$ is bijective.
\beqn
\label{Eq:fundamental-diagram}
\xymatrix@M+0.2em{
1 \ar[r] &
G(\hat W(R[F^n])^\oplus)_\mu \ar[r]^-{\alpha_\mu} \ar[d]^{\gamma_{\bar g}} &
G(\sN_n(R)^\oplus)_\mu \ar[r]^-{\pi_\mu} 
\ar[d]_\cong^{\gamma_{g}} &
G(\bar N_n(R)^\oplus)_\mu \ar[r] \ar[d]^{\bar\gamma_{g}} & 1 \\
1 \ar[r] &
G(\hat W(R[F^n])) \ar[r]^-\alpha &
G(N_n(R)^{\nil}) \ar[r]^-\pi &
G(\bar N_n(R)) \ar[r] & 1
}
\eeqn
Here $\hat W(R[F^n])^\oplus$ is considered as an ideal of $B_n(R)^\oplus$ as before.
\eCo

\bproof
The vertical maps are different instances of the orbit map \eqref{Eq:gamma-g}, thus $\bar\gamma_g(x)=\tau(x)^{-1}\sigma_{g}(x)$ with $\sigma_g=\Ad(g)\circ\sigma$ etc. Here $\gamma_g$ is well-defined as a function of $g$ and $\gamma_{\bar g}$ is well-defined as a function of $\bar g$ by \eqref{Eq:GammaGmuAK-functor-quot}, using that the ideal $W(R[F^n])^\oplus$ of $B_n(R)^\oplus$ is a $W_n(R)^\oplus$-module, and the ideal $\sN_n(R)^\oplus$ is an $A_n(R)^\oplus$-module. The diagram commutes since \eqref{Eq:gamma-g} is functorial with respect to the underlying frame and ideal.
The exact rows are given by Proposition \ref{Pr:rhon-surjective-full}. The map $\gamma_{g}$ is bijective by Lemma \ref{Le:dot-sigma-loc-nilpotent} since the ideal $ \sN_n(R)^\oplus$ is $\dot\sigma$-nilpotent.
\eproof

\bPr
\label{Pr:rhon-surjective-full}
The functor $\varrho_n^{G,\mu}$ of \eqref{Eq:rhon-Gamma} is surjective on objects and full.
\ePr

\bproof
The functor is surjective on objects by Lemma \ref{Le:Gpi0-Gpimu-surjective}.
By Lemma \ref{Le:functor-pi-G-mu} we have to show that for each $g\in G(A_n(R))$ the map $\bar \gamma_g$ in \eqref{Eq:fundamental-diagram} is surjective. This follows from the right square of the diagram.
\eproof

\bDe
Assume that $g\in G(A_n(R))$ with image $\bar g\in G(W_n(R))$ is given.
We define a map of pointed sets
\[
\psi_{g}=\pi_\mu\circ(\gamma_{g})^{-1}\circ\alpha:
G(\hat W(R[F^n]))\to G(\bar N_n(R)^\oplus)_\mu 
\]
as a composition in \eqref{Eq:fundamental-diagram}.
\eDe

\bLe
\label{Le:sequence-pointed-sets}
There is an exact sequence of pointed sets
\beqn
\label{Eq:sequence-pointed-sets}
1\to G(\hat W(R[F^n])^\oplus)_\mu
\xrightarrow{\;\gamma_{\bar g}\;}
G(\hat W(R[F^n]))
\xrightarrow{\;\psi_{g}\;}
E_n^{G,\mu}(R,g)\to 1.
\eeqn
\eLe

The group $E_n^{G,\mu}(R,g)$ was defined in \eqref{Eq:EGmunRg}. All objects in \eqref{Eq:sequence-pointed-sets} are groups, but a priori the arrows are only maps of pointed sets.

\bproof
Direct verification, using that $E_n^{G,\mu}(R,g)=(\bar\gamma_g)^{-1}(1)$.
\eproof

\bLe
\label{Le:left-column-abelian}
In the diagram \eqref{Eq:fundamental-diagram}, the groups in the left column are abelian, $\gamma_{\bar g}$ is a group homomorphism, and the left column is isomorphic to the complex of abelian groups 
$C^{\Fg,\mu}(\u{\hat W}(R[F^n]),\bar g)$ in degrees $\{-1,0\}$ defined in \eqref{Eq:ZngmuRg}.
\eLe

\bproof
The ideal $\hat W(R[F^n])^\oplus$ of $B_n(R)^\oplus$ has square zero by Proposition \ref{Pr:ker-tilde-rho-n-s}, so the groups in the left column are abelian, and $\gamma_{\bar g}$ is a difference of homomorphisms, thus a homomorphism. 
The inverse of the isomorphism $\exp$ of \eqref{Eq:exp-Gammagmu-GammaGmu} for $\u A=\u B_n(R)$ and $K=\hat W(R[F^n])^\oplus$
identifies the left column with the complex
$C^{\Fg,\mu}(\u{\hat W}(R[F^n]),\bar g)$
(using also Remark \ref{Re:CgmuNg-CgmuN'g'}).
\eproof

\bLe
\label{Le:psi-homomorphism}
In the diagram \eqref{Eq:fundamental-diagram}, the map $(\gamma_{g})^{-1}\circ\alpha$ is a group homomorphism, and hence $\psi_{g}$ is a group homomorphism.
\eLe

\bproof
It suffices to prove the first assertion because $\pi_\mu$ is a group homomorphism. We have $(\gamma_{g})^{-1}=\tau^{-1}\circ(\dot\gamma_{g})^{-1}$ where 
$\tau:G(\sN_n(R)^\oplus)_\mu\xrightarrow\sim G(N_n(R)^{\nil})$ is an isomorphism of groups, and $(\dot\gamma_{g})^{-1}$ is given by \eqref{Eq:dot-gamma-g-inverse}. By that formula it suffices to show that for $x,x'\in G(W(R[F^n]))$ and $m,m'\ge 0$ the elements 
\[
y=\dot\sigma_{g}^m(\alpha(x))
\quad\text{and}\quad
y'=\dot\sigma_{g}^{m'}(\alpha(x'))
\]
of $G(N_n(R)^{\nil})$ commute. Since $\dot\sigma_{g}$ is a group homomorphism we may assume that $m'=0$. Then $y'\in G(W(R[F^n]))$, and the commutator $[y,y']$ vanishes by Lemma \ref{Le:commutator} since the product of ideals $N_n(R)^{\nil}\cdot W(R[F^n])$ in $B_n(R)$ is zero by Proposition \ref{Pr:ker-tilde-rho-n}.
\eproof

\bLe
\label{Le:commutator}
Let $S$ be a\/ $\ZZ_p$-algebra and let $I,J\subseteq S$ be ideals. For $x\in G(I)$ and $y\in G(J)$ the commutator $[x,y]=xyx^{-1}y^{-1}$ lies in $G(IJ)$.
\eLe

\bproof
Let $G=\Spec C$ and $K\subseteq C$ the augmentation ideal. The commutator map $G\times G\to G$ sends $G\times\{e\}$ and $\{e\}\times G$ to $\{e\}$, hence the associated homomorphism $C\to C\otimes C$ restricts to $K\to K\otimes K$.
\eproof

\bproof[Proof of Proposition \ref{Pr:Main}]

The functor $\varrho_{n}^{G,\mu}$ is surjective on objects and full by Proposition \ref{Pr:rhon-surjective-full}.
Lemmas \ref{Le:sequence-pointed-sets}, \ref{Le:left-column-abelian}, and \ref{Le:psi-homomorphism} give an isomorphism $E_n^{G,\mu}(R,g)\cong H^0(C^{\Fg,\mu}(\u{\hat W}(R[F^n]),\bar g))$. Clearly this isomorphism is functorial in $R$. It is functorial in $g$ because the entire diagram 
\eqref{Eq:fundamental-diagram} for varying $g\in G(A_n(R))$ becomes a functor on $\pWin^{G,\mu}(\u A_n(R))$ as follows. To a morphism $g*h\to h$ in this groupoid defined by $h\in G(A_n(R)^\oplus)_\mu$  we associate the morphism from \eqref{Eq:fundamental-diagram} for $g*h$ to \eqref{Eq:fundamental-diagram} for $g$ defined by $\Ad'(\bar h)$, $\Ad'(h)$, $\Ad(h)$ in the upper line and by $\Ad'(\tau(\bar h))$, $\Ad'(\tau(h))$, $\Ad(\tau(h))$ in the lower line; here $\Ad'$ is provided by Lemma \ref{Le:Ad'}. 
\eproof


\section{Formal groups and $n$-smooth group schemes}
\label{Se:Zink-functor}

We recall the central construction of \cite{Zink:The-Display} with generalisations by \cite{Lau-Zink} and \cite{Drinfeld:The-Lau}, using the frame formalism employed here. Let $R$ be an $\FF_p$-algebra.

\subsection{Locally nilpotent algebras}

We denote by $\ttNil_R$ the category of locally nilpotent commutative $R$-algebras (non-unitary) and by $\ttNil_{n,R}$ be the category of all $N\in\ttNil_R$ such that $x^{p^n}=0$ for all $x\in N$.
For $N\in\ttNil_R$ we form the augmented $R$-algebra $R\oplus N$ with multiplication $(a\oplus x)(b\oplus y)=ab\oplus(ay+bx+xy)$ for $a,b\in R$ and $x,y\in N$. 

\subsubsection*{Modules associated to locally nilpotent algebras}

For $N\in\ttNil_R$ the frame $\u W(R\oplus N)$ 
contains the ideals $\hat W(N)^\oplus\subseteq W(N)^\oplus$; see \S\ref{Se:Witt-crys-frames.ideals-Witt}. The associated $\u W(R\oplus N)$-modules will be denoted by $\u{\hat W}(N)\subseteq\u W(N)$. They become modules over $\u W(R)$ by restriction of scalars.
If $N$ lies in $\ttNil_{n,R}$ then $\u{\hat W}(N)$ and $\u W(N)$ are modules over $\u {W\!}_n(R)$; see Remark \ref{Re:Wb-module-WnR}.

\subsection{Formal and $n$-smooth groups}

Let $\LG(R)=\LG_{\infty}(R)$ denote the category of commutative formal Lie groups over $R$, and for $n<\infty$ let $\LG_n(R)$ denote the category of commutative $n$-smooth group schemes over $R$ as defined in \cite[Definition 2.2.2]{Drinfeld:The-Lau}. For $0\le n\le\infty$ we have a fully faithful additive functor
\beqn
\label{Eq:LGn-FunNilnAb}
\LG_n(R)\to\Fun(\ttNil_{n,R},\Ab), \quad G\mapsto G_{\nil},
\eeqn
\beqn
G_{\nil}(N)=\ker(G(R\oplus N)\to G(R)).
\eeqn
Here one could also use the category $\Nil_{n,R}$ of nilpotent algebras instead of $\ttNil_{n,R}$; this makes no difference since the functors $G_{\nil}$ are limit-preserving.
For $m\le n$ the truncation functor $\LG_n(R)\to\LG_m(R)$, $G\mapsto G[F^m]$ corresponds to the restriction of functors along the inclusion $\ttNil_{m,R}\subseteq\ttNil_{n,R}$. One recovers $G$ from $G_{\nil}$ as follows: For an $R$-algebra $S$ we have
\beqn
\label{Eq:G-Gnil}
G(S)=\begin{cases}G_{\nil}(S[F^n])&\text{if $n<\infty$} \\
\bigcup_{m<\infty}G_{\nil}(S[F^m])&\text{if $n=\infty$}. \end{cases}
\eeqn

\subsection{The Zink complex and Zink functor}

\bDe
We define a functor
\beqn
Z_R:\Disp(R)\to\Fun(\ttNil_R,Ch(\Ab)),
\eeqn
\beqn
Z_R(\u M)(N)=\Gamma_1(\u M\otimes\u{\hat W}(N)),
\eeqn
using the complex $\eqref{Eq:Ci}$ with $i=1$. We will call $Z_R(\u M)(N)$ the Zink complex.
\eDe

Explicitly
\beqn
Z_R(\u M)(N)=[(M\otimes_{W(R)^\oplus}\hat W(N)^\oplus)_1\xrightarrow{\;\gamma_F\;}M^\tau\otimes_{W(R)}\hat W(N)]
\eeqn
with $\gamma_F(x)=F(x^\sigma)-x^\tau$.
If $\u M$ is effective, following \eqref{Eq:Ci-eff} we have
\beqn
\label{Eq:ZRMN-effective}
Z_R(\u M)(N)=[(M\otimes_{W(R)^\oplus}\hat W(N)^\oplus)_1\xrightarrow{\;\phi_1-t\;}M_0\otimes_{W(R)}\hat W(N)]
\eeqn
where $\phi:M\to M_0$ is the $\sigma$-linear map corresponding to $F$.

\bRe
\label{Re:Zink-complex-classical}
If the graded module $M$ is generated by degrees $0,1$, then $\u M$  corresponds to the classical display $\PPP=(M_0,M_1,\phi_0,\phi_1)$, and $Z_R(\u M)$ is the complex $[G_{\PPP}^{-1}\to G_{\PPP}^0]$ of \cite{Zink:The-Display}.
\eRe

\bDe
We define a functor
\beqn
Z_{n,R}:\Disp_n(R)\to\Fun(\ttNil_{n,R},Ch(\Ab)),
\eeqn
\beqn
\label{Eq:ZnRMN}
Z_{n,R}(\u M)(N)=\Gamma_1(\u M\otimes\u{\hat W}(N)).
\eeqn
\eDe

Explicitly
\[
Z_{n,R}(\u M)(N)=[(M\otimes_{W_n(R)^\oplus}\hat W(N)^\oplus)_1\xrightarrow{\gamma_F}M^\tau\otimes_{W_n(R)}\hat W(N)]
\]
with $\gamma_F$ as above. If $\u M$ is effective, the obvious analogue of \eqref{Eq:ZRMN-effective} holds.

\bPr
\label{Pr:Zink-functor}
Let $\u M$ be an effective display over $R$. Then for every $N\in\ttNil_R$ we have $H^{-1}(Z_R(\u M)(N))=0$, and the functor 
\[
H^0\circ Z_R(\u M)\in\Fun(\ttNil_R,\Ab)
\] 
lies in the image of $\LG(R)$ under the functor \eqref{Eq:LGn-FunNilnAb} for $n=\infty$.
\ePr

\bproof
If $M$ is generated by degrees $0,1$, this is a reformulation of \cite[Theorem 81]{Zink:The-Display}; see Remark \ref{Re:Zink-complex-classical}. As pointed out in \cite[\S 7.1.3]{Drinfeld:The-Lau}, the proof works in greater generality. Let us indicate the key points, using the present formalism. Clearly the functor $Z_R(\u M):\ttNil_{R}\to Ch(\Ab)$ is exact and limit preserving. The category $\Mod_R$ of $R$-modules can be identified with the category of all $N\in\Nil_R$ with $N^2=0$. Following \cite{Zink:The-Display} it suffices to show that $H^{-1}(Z_R(\u M)(N))=0$ when $N^2=0$ and that the functor $\Mod_R\to\Ab$, $N\mapsto H^0(Z_R(\u M)(N))$ is represented by the finite projective $R$-module $V=M_0/tM_1$. If $N^2=0$ then $\hat W(N)=\bigoplus_{i\ge 0}N$ because all higher terms in the Witt vector addition vanish. On can define an inclusion of $\u W(R)$-modules $\u{\hat W}(N)\to\u{\hat W}(N)'=\u U$ which is the identity in degrees $\le 0$ and the inclusion of $W(R)$-modules $\hat I(R)\to\hat W(R)$ in degrees $>0$. The multiplication $t:U_1\to U_0$ is the identity, the multiplication $t:U_i\to U_{i-1}$ is zero for $i>1$, and $\sigma_i:U_i\to U_0$ is the shift map $\dot\sigma:\hat W(N)\to\hat W(N)$, $(a_0,a_1,\ldots)\mapsto(a_1,a_2,\ldots)$. Now we have an exact sequence
\[
0\to Z_R(\u M)(N)\to \Gamma_1(\u M\otimes\u{\hat W}(N)')\to V\otimes_RN[1]\to 0
\]
where $V$ is the finite projective $R$-module $M_0/tM_1$. The complex in the middle is acyclic since it takes the form $[X\xrightarrow{1-u}X]$ with $X=M_0\otimes_{W(R)}\hat W(N)$ where $u=\phi\otimes\dot\sigma$ is pointwise nilpotent since $\dot\sigma$ is pointwise nilpotent on $\hat W(N)$.
\eproof

Proposition \ref{Pr:Zink-functor} gives an additive functor
\beqn
\FZ_R:\Disp^{\eff}(R)\to\LG(R),\quad\FZ_R(\u M)_{\nil}=H^0\circ Z_R(\u M).
\eeqn
The construction works when $R$ is a ring in which $p$ is nilpotent. This is Zink's functor $BT$ of \cite{Zink:The-Display}; the notation $\FZ_R$ was suggested by Drinfeld in \cite{Drinfeld:The-Lau}.

\bCo
\label{Co:Zink-functor-n}
Let $\u M$ be an effective $n$-truncated display over $R$. Then for every $N\in\ttNil_{n,R}$ we have $H^{-1}(Z_{n,R}(\u M)(N))=0$, and the functor 
\beqn
H^0\circ Z_{n,R}(\u M)\in\Fun(\ttNil_{n,R},\Ab)
\eeqn
lies in the image of $\LG_n(R)$ under the functor \eqref{Eq:LGn-FunNilnAb}.
\eCo

\bproof
Following  \cite[Proposition 3.11]{Lau-Zink} this is a direct consequence of Proposition \ref{Pr:Zink-functor}: The truncated display $\u M$ lifts to a display $\u N$, and the functor $Z_{n,R}(\u M)$ is the restriction of $Z_{n,R}(\u N)$ under $\ttNil_{n,R}\to\ttNil_R$.
\eproof

Corollary \ref{Co:Zink-functor-n} gives an additive functor
\beqn
\label{Eq:FZnR-LGnR}
\FZ_{n,R}:\Disp_n^{\eff}(R)\to\LG_n(R),\quad\FZ_{n,R}(M,F)_{\nil}=H^0\circ Z_{n,R}(\u M).
\eeqn

\subsection{Cohomology of $n$-smooth group schemes}

Let $H\in\LG_n(R)$ for an $\FF_p$-algebra $R$.

\bLe
\label{Le:et-coh-n-smooth=0}
We have $H^i_{\et}(\Spec R,H)=0$ for $i>0$.
\eLe

\bproof
Using the exact sequences $0\to H[F]\to H\to F^*(H[F^{n-1}])\to 0$ we may assume that $n=1$. 
It suffices to show that for every etale covering $R\to R'$ the \v Cech cohomology $H^i(R'/R,H)$ is trivial, i.e.\ the augmented \v Cech complex $\tilde C^*(R'/R,H)$ is acyclic. 
By [SGA 3, Expos\'e VII$_{\mathrm A}$] (see also \cite[\S2]{deJong:Finite-locally-free}) the group scheme $H$ corresponds to a pair $(M,\phi)$ where $M$ is a finite projective $R$-module and $\phi:M\to M$ is $F$-linear. By adding $(M',0)$ we can assume that $M$ is free. Then the Cartier dual $H^\vee$ is given by an exact sequence of fppf sheaves
\beqn
0\to H^\vee\to \GG_a^m\xrightarrow{\;F-A\;}\GG_a^m\to 0
\eeqn
for a matrix $A\in M_m(R)$. The functor $\u\Hom(-,\GG_m)$ gives an exact sequence of fppf sheaves
\beqn
\label{Eq:et-coh-n-smooth=0-resolution}
0\to(\hat W[F])^m\xrightarrow{\;V-A^t\;}(\hat W[F])^m\to H\to 0.
\eeqn
Indeed, exactness on the left is clear, and the the presheaf cokernel of $V-A^t$ is a $1$-smooth group scheme $H'$ (for example $V-A^t$ can be viewed as the Zink complex of a truncated display of level $1$). Then $H/H'$ is a $1$-smooth group scheme again. Since $\GG_a^{\vee\vee}=\GG_a$ it follows that $H^\vee\to H'^\vee$ is injective, hence $H=H'$. The argument shows that \eqref{Eq:et-coh-n-smooth=0-resolution} is exact as presheaves. Hence it suffices to show that the augmented \v Cech complex $\tilde C^*(R'/R,\hat W[F])$ is acyclic.
If $R\to R'$ is etale, then $W_n(R)\to W_n(R')$ is etale, and 
\[
\hat W(R'[F^n])=\hat W(R[F^n])\otimes_{W_n(R)}W_n(R').
\]
Indeed, for every $N\in\ttNil_{n,R}$ we have
\[
\hat W(N\otimes_RR')=\hat W(N)\otimes_{W_n(R)}W_n(R')
\]
because both sides are limit preserving, so we can assume that $N$ is nilpotent, even $N^2=0$ by induction; then $\hat W(N)=\bigoplus_{i\ge 0}N$. Hence the Amitsur complex gives the result.
\eproof

\bLe
\label{Le:Hifppf-Hisyn-nsmooth}
If $R$ is semiperfect then $H^i_{\fppf}(\Spec R,H)=H^i_{\syn}(\Spec R,H)=0$ for $i>0$.
\eLe

\bproof
For $\tau=\fppf$ or $\tau=\syn$ let $\pi:R_{\tau}\to R_{\et}$ be the morphism of sites such that $\pi_*$ is the restriction of sheaves to the small etale site. By Lemma \ref{Le:et-coh-n-smooth=0} it suffices to show that $R^i\pi_*H=0$ for $i>0$. As in the proof of Lemma \ref{Le:et-coh-n-smooth=0} we can assume that $n=1$ and that $H$ corresponds to a pair $(M,F)$ where $M$ is free. 
Since $R$ is semiperfect there is a pair $(M',F')$ with $F^*(M',F')=(M,F)$ and hence $H=F^*(H')$ for an $H'\in\LG_1(R)$. After Zariski localisation we find an exact sequence $0\to H'\to A'\to B'\to 0$ where $A$, $B$ are commutative smooth group schemes (Raynaud), which gives under $F^*$ a similar sequence $0\to H\to A\to B\to 0$. By \cite[\S11]{Grothendieck:BrauerIII} we have $R^i\pi_*A=R^i\pi_*B=0$ for $i>0$; the case $\tau=\syn$ is analogous to the case $\tau=\fppf$ treated there. Hence it suffices to show that $\pi_*A\to\pi_*B$ is surjective. Since $H'$ is annihilated by $F$, the relative Frobenius $A'\to A$ factors as $A'\to B'\to A$, and the composition $B'\to A\to B$ is the relative Frobenius of $B$. The homomorphism $B'(R)\to B(R)$ can be identified with the endomorphism of $B(R)$ induced by $F:R\to R$.
Since $R$ is semiperfect and $B'$ is smooth, this endomorphism is surjective. Hence $A(R)\to B(R)$ is surjective, which remains true after etale extension of $R$, and thus $\pi_*A\to\pi_*B$ is even surjective as presheaves.
\eproof

\section{Drinfeld's group scheme}
\label{Se:Drinfeld-group-scheme}

In this section we recall a central construction of \cite{Drinfeld:The-Lau} and relate it to Proposition \ref{Pr:Main}.
Let $(G,\mu)$ be a pair as in \S\ref{Se:displays-Gmu-groupoid} where $\mu$ is $1$-bounded. Let $\tilde\mu={\Ad}\circ\mu$ as earlier.
The functor $\Lie^{G,\mu}$ of \eqref{Eq:tildeAdGmu} induces a homomorphism of fibered categories over $\Aff_k$
\[
\Lie^{G,\mu}:
\Disp_n^{G,\mu}\to 
\Disp_n\otimes_{\FF_p}k.
\]

\bLe
\label{Le:1bd-eff}
The shifted version $\Lie^{G,\mu}(-1)=(-1)\circ\Lie^{G,\mu}$ gives a functor
\[
\Lie^{G,\mu}(-1):
\Disp_n^{G,\mu}\to 
\Disp_n^{\eff}\otimes_{\FF_p}k.
\]
\eLe

\bproof
In the setting of \S\ref{Se:variations.adjoint} we have $\Fg_{\tilde\mu,i}=0$ for $i<-1$ since $\mu$ is $1$-bounded and since $\Fg_{\tilde\mu,i}=\Fg^{-i}$. Hence $\Fg_{\tilde\mu}(-1)$ is effective.
\eproof

The functor $\Lie^{G,\mu}(-1)$ is the functor (C.10) of \cite{Drinfeld:The-Lau}. It allows to define a functor
\beqn
\FD_n^{G,\mu}:\Disp_n^{G,\mu}\to\LG_n\otimes_{\FF_p}k,\qquad \FD_n^{G,\mu}=\FZ_n\circ\Lie^{G,\mu}.
\eeqn
This functor corresponds to a commutative $n$-smooth group scheme over $\Disp_n^{G,\mu}$, which will also be denoted by $\FD_n^{G,\mu}$. In \cite{Drinfeld:The-Lau} this group scheme is denoted $\Lau_n^{G,\mu}$. In formulas, if $Q$ denotes the universal $(G,\mu)$-display over $\Disp_n^{G,\mu}$, then
\beqn
\label{Eq:DnGmu}
\FD_n^{G,\mu}=\FZ_n(\Lie^{G,\mu}(Q)(-1)).
\eeqn

\bDe
Let $R$ be a $k$-algebra.
Each $g\in G(W_n(R))$ gives $g:\Spec R\to\Disp_n^{G,\mu}(R)$ and thus an $n$-smooth group scheme $g^*(\FD_n^{G,\mu})$ over $R$. We denote its global sections by
\beqn
D_n^{G,\mu}(R,g)=g^*(\FD_n^{G,\mu})(\Spec R).
\eeqn
This abelian group is functorial in $R$ and $g$.
\eDe

\bLe
\label{Le:DnGmuRg-H0ZngmuRg}
There is an isomorphism of abelian groups
$D_n^{G,\mu}(R,g)\cong H^0(C^{\Fg,\mu}(\u{\hat W}(R[F^n]),g)$
which is functorial in $R$ and $g$.
\eLe

The complex of abelian groups $C^{\Fg,\mu}(\u{\hat W}(R[F^n]),g)$ was defined in \eqref{Eq:ZngmuRg}.

\bproof
By tracing the definitions. With $\Fb_n=R[F^n]$ we have
\[
D_n^{G,\mu}(R,g)=\FZ_n(\Lie^{G,\mu}(g)(-1))_{\nil}(\Fb_n)=H^0(Z_{n,R}(\Lie^{G,\mu}(g)(-1))(\Fb_n))
\]
using \eqref{Eq:DnGmu}, \eqref{Eq:G-Gnil}, and \eqref{Eq:FZnR-LGnR}, moreover
\[
Z_{n,R}(\Lie^{G,\mu}(g)(-1))(\Fb_n)=\Gamma_1(\Lie^{G,\mu}(g)(-1)\otimes\u{\hat W}(\Fb_n))
\] 
\[
=\Gamma_0(\tilde\Ad{}^{G,\mu}(g)\otimes\u{\hat W}(\Fb_n))
=C^{\Fg,\mu}(\u{\hat W}(\Fb_n),g)
\]
where the equalities are \eqref{Eq:ZnRMN}, evident, and \eqref{Eq:ZngmuRg}. All isomorphisms are functorial.
\eproof

Now we have the following more useful interpretation of the isomorphism \eqref{Eq:Pr-Main}.

\bCo
\label{Co:Main}
For each $R\in\Aff_k^{\qrsp}$ and $g\in G(A_n(R))$ with image $\bar g\in G(W_n(R))$ there is an isomorphism of abelian groups
$E_n^{G,\mu}(R,g)\cong D_n^{G,\mu}(R,\bar g)$
which is functorial in $R$ and $g$. \qed
\eCo

\bproof[Proof of Theorem \ref{Th:Main-pre}]
Proposition \ref{Pr:Main} and Corollary \ref{Co:Main}.
\eproof


\section{Reconstruction of algebraic stacks}
\label{Se:Reconstruction}

\subsection{Topologies}
\label{Se:Reconstruction.topologies}

Let $k$ be a quasisyntomic $\FF_p$-algebra. The category of quasiregular semiperfect affine $k$-schemes and the category of quasisyntomic affine $k$-schemes and the inclusion functors will be denoted by
\[
\Aff_k^{\qrsp}\xrightarrow{\;v\;}\Aff^{\qsyn}_k\xrightarrow{\;u\;}\Aff_k,
\qquad j=uv.
\]
On these categories we consider the etale, ind-syntomic, and quasisyntomic topology. In each case, coverings are finite families of morphisms of the given type which are jointly surjective; a ring homomorphism $R\to R'$ will be called ind-syntomic if $R'$ is a filterd colimit of syntomic $R$-algebras. Ind-syntomic implies quasisyntomic. On $\Aff^{\qsyn}_k$ and $\Aff_k$ we also consider the syntomic topology, and on $\Aff_k$ also the fppf and fpqc topology. 

The functor $v$ induces an equivalence on $\tau$-sheaves and on $\tau$-stacks for $\tau\in\{\indsyn,\qsyn\}$ because for every $k$-algebra $R$ there is an ind-syntomic cover $R\to R'$ where $R'$ is semiperfect.

Algebraic stacks over $k$ will be viewed as fibered categories over $\Aff_k$.

\bRe
\label{Re:algebraic-stacks}
If $\XXX$ is an algebraic space over $k$ or an algebraic stack over $k$ with quasi-affine diagonal, then $\XXX$ is an fpqc stack and hence a quasisyntomic stack by \cite[\href{https://stacks.math.columbia.edu/tag/0APL}{Proposition 0APL}, \href{https://stacks.math.columbia.edu/tag/0GRH}{Proposition 0GRH}]{Stacks} or \cite[Corollaire 10.7]{LMB}. If $\XXX$ is an algebraic stack locally of finite presentation over $k$, then $\XXX$ is an fppf stack by  \cite[Corollaire 10.7]{LMB} again, hence a syntomic stack, and $\XXX$ is limit preserving by \cite[Proposition 4.15(i)]{LMB}, hence an ind-syntomic stack.
\eRe

\subsection{Reconstruction of syntomic algebraic stacks}

For a presheaf or fibered category $\FFF$ over $\Aff_k$ we denote by $\res(\FFF)=u^{-1}(\FFF)$ its restriction to $\Aff^{\qsyn}_k$.
The following is a variant of \cite[Proposition A.0.1]{Drinfeld:On-Shimurian} with the same proof.

\bPr
\label{Pr:res-fully-faithful}
Let $\XXX$ be a syntomic algebraic stack over $k$ and let $\FFF$ be an etale stack over $\Aff_k$.\footnote{Following standard terminology, a syntomic algebraic stack is an algebraic stack which is syntomic, and an etale stack is a stack for the etale topology.}  Then the restriction functor
\begin{equation}
\label{Eq:res-XXX-FFF}
\res_{\XXX,\FFF}:\Fun_{\Aff_k}(\XXX,\FFF)\to\Fun_{\Aff^{\qsyn}_k}(\res(\XXX),\res(\FFF))
\end{equation}
is an equivalence of categories. \qed
\ePr

\bDe
\label{De:stack-algebraic}
An etale stack $\XXX'$ over $\Aff^{\qsyn}_k$ will be called syntomic algebraic if $\XXX'=\res(\XXX)$ for a syntomic algebraic stack $\XXX$ over $k$, which is essentially unique by Proposition \ref{Pr:res-fully-faithful}. We will say that $\XXX'$ has a certain property of algebraic stacks over $k$ if this holds for $\XXX$. Similarly, a morphism between syntomic algebraic stacks over $\Aff^{\qsyn}_k$ will be said to have a certain property of morphisms of algebraic stacks if its unique extension to a morphism of syntomic algebraic stacks over $\Aff_k$ has this property.
\eDe

\bRe
If an etale stack $\XXX'$ over $\Aff^{\qsyn}_k$ is syntomic algebraic, then $\XXX'$ is an ind-syntomic stack, and even a quasisyntomic stack if under additional conditions as in Remark \ref{Re:algebraic-stacks}.
\eRe

\subsection{Reconstruction of the gerbe condition}

\subsubsection*{Classifying stacks}

For a sheaf of groups $H$ on a site $(\SSS,\tau)$ we denote by $B_\tau H$ the stack $H$-torsors with respect to $\tau$. For example, if $G$ is a flat group scheme locally of finite presentation over $k$ one defines $BG=B_{\fppf}G$ over $\Aff_k$; this is a smooth algebraic stack over $k$ by \cite[\href{https://stacks.math.columbia.edu/tag/0DLS}{Lemma 0DLS}]{Stacks}. 

\bLe
\label{Le:res-BL-gen}
Let $G$ be an affine flat group scheme of finite presentation over $k$, and let $\tau\in\{\syn,\indsyn,\qsyn\}$. Then $\res(BG)\cong B_{\tau}\res(G)$ as fibered categories over $\Aff^{\qsyn}_k$, in particular $B_{\tau}\res(G)$ is a syntomic algebraic stack and smooth (Definition \ref{De:stack-algebraic}).
\eLe

\bproof
We have $BG=B_{\tau}G$ since $G$ is syntomic over $k$ by \cite[Proposition 27.26]{GW2}.
\eproof

\bPr
\label{Pr:res-gerbe}
Let $f:\XXX\to\YYY$ be a morphism of syntomic algebraic stacks over $k$ and 
$G$ a commutative affine flat group scheme of finite presentation over $\YYY$. The following are equivalent.
\begin{enumerate}
\item
\label{It:f-fppf-gerbe}
The morphism $f$ is an fppf gerbe banded by $G$,
\item
\label{It:res(f)-syn-gerbe} 
the morphism $\res(f)$ is a syntomic gerbe banded by $\res(G)$,
\item
\label{It:res(f)-indsyn-gerbe}
the morphism $\res(f)$ is a quasisyntomic gerbe banded by $\res(G)$.
\end{enumerate}
\ePr

\bproof
If $f$ is an fppf gerbe banded by $G$, then $f$ has sections etale locally, and two sections are isomorphic syntomic locally since $G$ is syntomic, so $f$ is a syntomic gerbe banded by $G$. This proves \eqref{It:f-fppf-gerbe}$\Rightarrow$\eqref{It:res(f)-syn-gerbe}, and \eqref{It:res(f)-syn-gerbe}$\Rightarrow$\eqref{It:res(f)-indsyn-gerbe} is clear.
Assume that \eqref{It:res(f)-indsyn-gerbe} holds.
Let $(\Spec R_i\to\YYY)_i$ be a syntomic cover, thus $\Spec R_i\in\Aff^{\qsyn}_k$. There is a quasisyntomic cover $(\Spec k_{ij}\to\Spec R_i)_j$ in $\Aff^{\qsyn}_k$ such that the fibered category $\res(\XXX_{k_{ij}})$ over $\Aff^{\qsyn}_{k_{ij}}$ is a neutral quasisyntomic gerbe, thus isomorphic to $B_{\indsyn}\res(G_{k_{ij}})\cong\res(BG_{k_{ij}})$ using Lemma \ref{Le:res-BL-gen}, and then $\XXX_{k_{ij}}\cong BG_{k_{ij}}$ by Proposition \ref{Pr:res-fully-faithful}. We claim that $\pi:\XXX\to\YYY$ is smooth and surjective and the diagonal $\Delta:\XXX\to\XXX\times_\YYY\XXX$ is affine, syntomic, and surjective. Indeed, these statements hold after base change under $k\to k_{ij}$ and satisfy fpqc descent. Hence $\XXX$ is a syntomic gerbe over $k$. Moreover the inertia $\III_{\XXX/\YYY}$ and $f^*(G)$ are affine syntomic group schemes over $\XXX$ and thus syntomic over $k$; hence the given isomorphism $\res(\III_{\XXX/\YYY})\cong\res(f^*(G))$ of groups over $\res(\XXX)$ comes from a unique isomorphism $\III_{\XXX/k}\cong f^*(G)$ of group schemes over $\XXX$ by Proposition \ref{Pr:res-fully-faithful} again.
\eproof

\subsection{Limit preserving stacks}

We need a technical preparation for \S\ref{Se:Reconstruction.descent-algebraicity}.

\bDe
\label{De:limit-preserving}
A fibered category $\FFF$ over $\Aff^{\qsyn}_k$ will be called limit preserving if for each filtered colimit of $k$-algebras $R^{(\infty)}=\colim_j R^{(j)}$ such that $\Spec R^{(j)}\in\Aff^{\qsyn}_k$ for all $j$ we have 
\[
\FFF(R^{(\infty)})\cong\colim_j\FFF(R^{(j)}).
\]
\eDe

\bLe
\label{Le:limit-preserving-gen}
Let $\FFF$ be a $\tau$-stack over $\Aff^{\qsyn}_k$ with
$\tau\in\{\et, \syn, \indsyn, \qsyn\}$.
If for some $\tau$-cover $(\Spec k_i\to\Spec k)_i$ the restrictions $\FFF_{k_i}$ over $\Aff^{\qsyn}_{k_i}$ are limit preserving, then so is $\FFF$. 
\eLe

\bproof
Let $R^{(\infty)}=\colim_j R^{(j)}$ as in Definition \ref{De:limit-preserving}. The family
$(\Spec R^{(j)}\otimes_kk_i\to\Spec R^{(j)})_i$ is a $\tau$-cover for each $j$, including $j=\infty$. By descent it suffices to show that for each $k$-algebra of the form $k'=k_{i_1}\otimes_k\ldots\otimes_kk_{i_r}$ with $r\ge 1$, the natural functor
\[
\colim_j\FFF(R^{(j)}\otimes_kk')\to\FFF(R^{(\infty)}\otimes_kk')
\]
is an equivalence. This holds because the homomorphisms of $k$-algebras $k_{i_1}\to k'\to R^{(j)}\otimes_kk'$ give equivalences 
\[
\FFF(R^{(j)}\otimes_kk')\cong\FFF_{k_{i_1}}(R^{(j)}\otimes_kk')
\]
which are compatible with change of $j$.
\eproof

\subsection{Localisation}
\label{Se:Reconstruction.localisation}

This is another technical preparation for \S\ref{Se:Reconstruction.descent-algebraicity}.
Let $k\to R$ be a homomorphism of quasisyntomic $\FF_p$-algebras. The categories $\Aff_{R}$ and $\Aff^{\qsyn}_{R}$ can be identified with the relative categories $(\Aff_k/\Spec R)$ and $(\Aff^{\qsyn}_k/\Spec R)$, and there is the following commutative diagram with horizontal equivalences as in \cite[\href{https://stacks.math.columbia.edu/tag/0791}{Lemma 0791}]{Stacks}; here $\Spec R$ is viewed as a sheaf on $\Aff_k$, and $\tau$ can be discrete, etale, syntomic, ind-syntomic, or quasisyntomic.
\beqn
\label{Eq:Localisation}
\xymatrix@M+0.2em{
Sh_\tau(\Aff_{R}) \ar[r]^-\lambda_-\sim \ar[d]^{\res} &
Sh_\tau(\Aff_k)/\Spec R \ar[d]^{\res} \\
Sh_\tau(\Aff^{\qsyn}_{R}) \ar[r]^-{\lambda'}_-\sim &
Sh_\tau(\Aff^{\qsyn}_k)/\res(\Spec R)
}
\eeqn
Explicitly, if $i^{-1}:Sh_{\tau}(\Aff_k)\to Sh_{\tau}(\Aff_{R})$ is the restriction functor, the inverse of $\lambda$ sends $(\FFF\to\Spec R)$ to the fiber product of the diagram $i^{-1}(\FFF)\to i^{-1}(\Spec R)\leftarrow *$; the last map sends $(\Spec S\xrightarrow v\Spec R)\in\Aff_{R}$ to $v$. 
The diagram \eqref{Eq:Localisation} extends from $\tau$-sheaves to $\tau$-stacks.

\bLe
\label{Le:localisation-algebraic}
Let $k\to R$ by a syntomic ring homomorphism. Let $\XXX_R$ be an etale stack over $\Aff^{\qsyn}_{R}$ and let $\lambda(\XXX_R)=(\XXX\xrightarrow{\;u\;}\Spec R)$ be the corresponding etale stack over $\Aff^{\qsyn}_{k}$ lying over $\Spec R$. Then $\XXX_R$ is syntomic algebraic iff $\XXX$ is syntomic algebraic and $u$ is syntomic; see Definition \ref{De:stack-algebraic}.
\eLe

\bproof
Syntomic algebraic stacks over $R$ are equivalent to algebraic stacks $\XXX$ over $k$ together with a syntomic morphism $\XXX\to\Spec R$, and  $\XXX$ is necessarily syntomic over $k$. This equivalence is given by the stack version of the above functor $\lambda$; now use ${\res}\circ\lambda=\lambda'\circ{\res}$.
\eproof

\subsection{Descent of algebraicity}
\label{Se:Reconstruction.descent-algebraicity}

If $(\SSS,\tau)$ is a site,
for a presheaf of groupoids $X_\bullet$ on $\SSS$ we denote by $[X_\bullet]_\tau$ the associated $\tau$-stack.

\bLe
\label{Le:descend-algebraic}
Let $\XXX'\to\YYY'$ be a morphism of fibered categories over $\Aff^{\qsyn}_k$ where $\XXX'$ is a syntomic stack and $\YYY'$ is a syntomic algebraic stack, thus $\YYY'=\res(\YYY)$ for a syntomic algebraic stack $\YYY$ over $k$. Assume that there is a syntomic cover $\alpha:Y\to\YYY$ where $Y$ is a scheme, with restriction $\alpha':Y'\to\YYY'$ over $\Aff^{\qsyn}_k$, such that the fibered category $\ZZZ'=\XXX'\times_{\YYY'}Y'$ over $\Aff^{\qsyn}_k$ is a syntomic algebraic stack. Then $\XXX'$ is a syntomic algebraic stack.
\eLe

\bproof
The assumption means that $\ZZZ'=\res(\ZZZ)$ for a syntomic algebraic stack $\ZZZ$ over $k$. Let $\beta:X_0\to\ZZZ$ be a smooth cover where $X_0$ is a scheme and let $\beta':X_0'\to\ZZZ'$ be the restriction over $\Aff^{\qsyn}_k$. 
\[
\xymatrix@M+0.2em{
X'_0 \ar[r]^-{\beta'} & \ZZZ' \ar[r] \ar[d] & \XXX' \ar[d] \\
&  Y' \ar[r]^{\alpha'} & \YYY'
}
\]
Then $X_0'\to\XXX'$ is a syntomic epimorphism and thus $\XXX'\cong[X'_\bullet]_{\syn}$ where $X'_\bullet$ is the associated \v Cech groupoid given by $X'_1=X'_0\times_{\XXX'} X'_0$, etc. Now $X'_1$ is a syntomic algebraic space, namely $X'_1=\res(X_1)$ where
\[
X_1=X_0\times_{\ZZZ}(Y\times_{\YYY} X_0).
\]
The projections $X_1\rightrightarrows X_0$ are syntomic since $\alpha$ is syntomic and $\beta$ is smooth, thus $X'_\bullet=\res(X_\bullet)$ for a syntomic groupoid of syntomic algebraic spaces $X_\bullet$.
Let $\XXX=[X_\bullet]_{\fppf}$
be the associated syntomic algebraic stack. 
Then $\XXX=[X_\bullet]_{\syn}$ and hence $\res(\XXX)=[\res(X_\bullet)]_{\syn}=\XXX$.
\eproof

\bPr
\label{Pr:gerbe-algebraic}
Let $\YYY$ be a syntomic algebraic stack over $k$ and $G$ a commutative affine flat group scheme of finite presentation over $\YYY$. 
Assume that $f:\XXX'\to\YYY'=\res(\YYY)$ is a morphism of ind-syntomic stacks over $\Aff^{\qsyn}_k$ which is an ind-syntomic gerbe banded by $G'=\res(G)$.
Then there is a unique fppf gerbe $\XXX\to\YYY$ banded by $G$ with an isomorphism $\res(\XXX)\cong\XXX'$ as gerbes over $\YYY'$ banded by $G'$. In particular, $\XXX'$ is syntomic algebraic.
\ePr

\bproof
Initial remark and notation. For a morphism $\Spec R\to\YYY$ where $R$ is quasi-syntomic over $k$, the fiber product $\XXX'\times_{\YYY'}\res(\Spec R)$ is a fibered category over $\Aff^{\qsyn}_k$ equipped with a morphism to $\res(\Spec R)$, which is equivalent to a fibered category $\XXX'_R$ over $\Aff_R'$; see \S\ref{Se:Reconstruction.localisation}.

Let $(\Spec R_i\to\YYY)_i$ be a smooth cover, thus $\Spec R_i$ is syntomic over $k$ and lies in $\Aff^{\qsyn}_k$. 
There is an ind-syntomic cover $(\Spec k_{ij}\to\Spec R_i)_j$ over which $\XXX'_{R_i}$ has sections and hence $\XXX'_{k_{ij}}\cong B_{\indsyn}\res(G_{k_{ij}})\cong \res(BG_{k_{ij}})$ using Lemma \ref{Le:res-BL-gen}. The algebraic stack $BG_{k_{ij}}$ over $k_{ij}$ is limit preserving, see for example \cite[\href{https://stacks.math.columbia.edu/tag/0CMX}{Lemma 0CMX}]{Stacks}. Hence $\XXX'_{k_{ij}}$ is limit preserving in the sense of Definition \ref{De:limit-preserving}, which carries over to $\XXX'_{R_i}$ by Lemma \ref{Le:limit-preserving-gen}.
Since the ring homomorphism $R_i\to k_{ij}$ is ind-syntomic, it follows that there is a syntomic cover $(\Spec A_{ij}\to\Spec R_i)_j$ over which $\XXX'_{R_i}$ has a section, 
and hence $\XXX'\times_{\YYY'}\res(\Spec A_{ij})$ is isomorphic to $\res(BG)\times_{\YYY'}\res(\Spec A_{ij})$ as a fibered category over $\Aff^{\qsyn}_k$.
Let $Y=\coprod_{i,j}\Spec A_{ij}$. 
Then $\ZZZ'=\XXX'\times_{\YYY'} \res(Y)$ is isomorphic to $\res(BG)\times_{\YYY'}\res(Y)=\res(BG\times_{\YYY}Y)$,
in particular $\ZZZ'$ is syntomic algebraic,
and Lemma \ref{Le:descend-algebraic} implies that $\XXX'=\res(\XXX)$ for a syntomic algebraic stack $\XXX$ over $k$. The given morphism $\XXX'\to\YYY'$ extends uniquely to a morphism $\XXX\to\YYY$ by Proposition \ref{Pr:res-fully-faithful}, and this morphism is a gerbe banded by $G$ by Proposition \ref{Pr:res-gerbe}.
\eproof

\section{Geometric conclusions}
\label{Se:Geom-appl}

\subsection{Recollection}

We use the notation of \S\ref{Se:Reconstruction.topologies}.
By \S\ref{Se:Witt-crys-Gmu-disp} we have $1$-categorical presheaves of groupoids $\pWin^{G,\mu}(\u A_n)$ on $\Aff^{\qrsp}_k$ and
$\pWin^{G,\mu}(\u{W\!}_n)$ on $\Aff_k$, and
the etale stack associated to $\pWin^{G,\mu}(\u{W\!}_n)$ is the algebraic stack $\Disp_n^{G,\mu}$. There is a homomorphism 
\beqn
\label{Eq:rhonGmu-recollection}
\varrho_n^{G,\mu}:\pWin^{G,\mu}(\u A_n)\to j^{-1}(\pWin^{G,\mu}(\u{W\!}_n))
\eeqn
which is surjective on objects and full by Proposition \ref{Pr:rhon-surjective-full}, and its inertia is identified with the restriction of the $n$-smooth group scheme $\FD_n^{G,\mu}$ on $\Disp_n^{G,\mu}$ by Corollary \ref{Co:Main}. 

\subsection{The associated stack}

Let $\Disp_n^{G,\mu,\qsyn}\to\Aff_k^{\qsyn}$ and $\Disp_n^{G,\mu,\qrsp}\to\Aff_k^{\qrsp}$ be the restrictions of $\Disp_n^{G,\mu}$, and let
\beqn
\label{Eq:tDisp''-Affk''}
\ttDisp^{G,\mu,\qrsp}_n\to\Aff^{\qrsp}_k
\eeqn
denote the etale stack associated to $\pWin^{G,\mu}(\u A_n)$.
The homomorphism $\varrho_n^{G,\mu}$ induces a functor
\beqn
\label{Eq:rhon-tDisp''-Disp}
\varrho_n^{G,\mu}:\ttDisp^{G,\mu,\qrsp}_n\to \Disp^{G,\mu,\qrsp}_n
\eeqn
of fibered categories over $\Aff_k^{\qrsp}$,
which is an etale gerbe banded by $\FD_n^{G,\mu,\qrsp}=j^{-1}(\FD_n^{G,\mu})$ as a consequence of Proposition \ref{Pr:rhon-surjective-full} and Corollary \ref{Co:Main}.

\bLe
\label{Le:tDisp-indsyntomic}
The etale stack $\ttDisp^{G,\mu,\qrsp}_n$ is an ind-syntomic stack.
\eLe

\bproof
Since $\pWin^{G,\mu}(\u A_n)$ is a quasisyntomic sheaf of groupoids, $\ttDisp^{G,\mu,\qrsp}_n$ is a quasisyntomic prestack, i.e.\ morphisms form sheaves. We have to show that ind-syntomic descent is effective, and this follows from Lemma \ref{Le:Hifppf-Hisyn-nsmooth}. In more detail, let $\ZZZ\to \Aff^{\qrsp}_k$ be the ind-syntomic stack associated to $\pWin^{G,\mu}(\u A_n)$. The functor $\varrho_n^{G,\mu}$ induces a functor $\varrho_n:\ZZZ\to \Disp_n^{G,\mu,\qrsp}$. We have to show that each object $x$ of $\ZZZ(R)$ lies in $\pWin^{G,\mu}(\u A_n(R))$ up to isomorphism etale locally in $\Spec R$. After etale localisation we can assume that the image $\bar x=\varrho_n(x)$ is isomorphic to an object $\bar g\in G(W_n(R))$. Let $g\in G(A_n(R))$ be an inverse image of $\bar g$. Since the functor \eqref{Eq:rhonGmu-recollection} is full, the sheaf of isomorphisms $x\cong g$ lifting the given isomorphism $\bar x\cong\bar g$ has a section over an ind-syntomic cover $R\to R'$, so $x$ corresponds to a descent datum for $g$ with values in $\FD=\bar g^*(\FD_n^{G,\mu})$. But $H^1(R'/R,\FD)$ vanishes as a consequence of Lemma \ref{Le:Hifppf-Hisyn-nsmooth}, using that $H^1(R'/R,\FD)$ is compatible with filtered colimits of $R'$.
\eproof

\bRe
In view of Lemma \ref{Le:tDisp-indsyntomic} we could also define $\ttDisp^{G,\mu,\qrsp}_n\to\Aff^{\qrsp}_k$ as the ind-syntomic stack associated to $\pWin^{G,\mu}(\u A_n)$. A priori, this gives an ind-syntomic gerbe over $\Disp^{G,\mu,\qrsp}_n$, which is an etale gerbe by Lemma \ref{Le:tDisp-indsyntomic}.
\eRe

Lemma \ref{Le:tDisp-indsyntomic} implies that \eqref{Eq:tDisp''-Affk''} extends uniquely to an ind-syntomic stack
\beqn
\label{Eq:tDisp'-Affk'}
\ttDisp^{G,\mu,\qsyn}_n\to\Aff^{\qsyn}_k,
\eeqn
and the morphism \eqref{Eq:rhon-tDisp''-Disp} extends uniquely to a a morphism
\beqn
\label{Eq:rhon-tDisp'-Disp}
\varrho_n^{G,\mu}:\ttDisp^{G,\mu,\qsyn}_n\to \Disp^{G,\mu,\qsyn}_n
\eeqn
over $\Aff^{\qsyn}_k$, which is an ind-syntomic gerbe banded by $\FD_n^{G,\mu,\qsyn}=u^{-1}(\FD_n^{G,\mu})$.
Now Proposition \ref{Pr:gerbe-algebraic} implies:

\bCo
\label{Co:gerbe-final}
The fibered category \eqref{Eq:tDisp'-Affk'} is the restriction of a unique syntomic algebraic stack
\beqn
\label{Eq:tDisp-Affk}
\ttDisp^{G,\mu}_n\to\Aff_k,
\eeqn
and the morphism \eqref{Eq:rhon-tDisp'-Disp} extends uniquely to a morphism
\beqn
\label{Eq:tDisp-Disp}
\varrho_n^{G,\mu}:\ttDisp^{G,\mu}_n\to \Disp^{G,\mu}_n
\eeqn
over $\Aff_k$, which is a gerbe banded by the group scheme $\FD_n^{G,\mu}$.
\qed
\eCo 

\bRe
Since $\Disp_n^{G,\mu}$ is a quasicompact algebraic stack with affine diagonal, i.e.\ the stack associated to a smooth groupoid of affine schemes, the same holds for $\ttDisp_n^{G,\mu}$, in particular this is an fpqc stack. It follows that $\ttDisp_n^{G,\mu,\qsyn}$ and $\ttDisp_n^{G,\mu,\qrsp}$ are quasisyntomic stacks, in particular $\ttDisp_n^{G,\mu,\qrsp}$ is also the quasisyntomic stack associated to the sheaf of groupoids $\pWin^{G,\mu}(\u A_n)$.
\eRe

\subsection{Relation with the Shimurian BT stack}

Following Drinfeld \cite{Drinfeld:On-Shimurian}, Gardner--Madapusi \cite{Gardner-Madapusi} define a derived $p$-adic formal prestack $\BT_n^{G,\mu}$ over $\Spf W(k)$ and prove  that this is a quasicompact smooth $0$-dimensional $p$-adic formal Artin stack with affine diagonal (\cite[Theorem D]{Gardner-Madapusi}), which confirms  \cite[Conjecture C.3.1]{Drinfeld:On-Shimurian} that was formulated in the case $k=\FF_p$. Let $\bbBT_n^{G,\mu}\to\Spec k$ be the special fiber of $\BT_n^{G,\mu,}$. By \cite[Remark 9.1.3]{Gardner-Madapusi} there is a morphism $\phi_n^{G,\mu}:\bbBT_n^{G,\mu}\to\Disp_n^{G,\mu}$.

\bPr
\label{Pr:bbBT-ttDisp}
There is an isomorphism $\bbBT_n^{G,\mu}\cong\ttDisp_n^{G,\mu}$ of algebraic stacks over $k$. Under this isomorphism, the morphism $\phi_n^{G,\mu}$ coincides with the morphism $\varrho_n^{G,\mu}$ of \eqref{Eq:tDisp-Disp}.
\ePr

\bproof
For $\tau\in\{\qsyn,\qrsp\}$ let $\bbBT_n^{G,\mu,\tau}$ denote the restriction of $\bbBT_n^{G,\mu}$ to the category $\Aff^{\tau}_k$.
Here $\bbBT_n^{G,\mu,\qrsp}$ is isomorphic to the etale stack associated to the sheaf of groupoids $\pWin_n^{G,\mu}(\u A_n)$ over $\Aff^{\qrsp}_k$ by \cite[Lemma 9.2.3, Remark 6.9.5, and Remark 5.5.5]{Gardner-Madapusi}. The isomorphism $\bbBT_n^{G,\mu,\qrsp}\cong\ttDisp_n^{G,\mu,\qrsp}$ extends to an isomorphism $\bbBT_n^{G,\mu,\qsyn}\cong\ttDisp_n^{G,\mu,\qsyn}$ over $\Aff^{\qsyn}_k$ because both sides are ind-syntomic stacks, and this isomorphism extends to an isomorphism $\bbBT_n^{G,\mu}\cong\ttDisp_n^{G,\mu}$ over $\Aff_k$ by Proposition \ref{Pr:res-fully-faithful}, using that both sides are smooth algebraic stacks. This proves the first assertion. By the construction in \cite[Remark 9.1.3 and Lemma 9.2.3]{Gardner-Madapusi}, the restriction of $\phi_n^{G,\mu}$ to the prestack $\pWin_n^{G,\mu}(\u A)$ over $\Aff_k^{\qrsp}$ coincides with the functor $\varrho_n^{G,\mu}$ induced by the frame homomorphisms $\varrho_n:\u A_n(R)\to\u{W\!}_n(R)$ for quasiregular semiperfect $R$. There is no ambiguity because this frame homomorphism is unique (determined either by the universal property of $A_{\crys}(R)$ or by the universal property of $W(R)$).
\eproof

\bCo[{\cite[Conjecture C.5.3]{Drinfeld:The-Lau}}]
The morphism $\phi_n^{G,\mu}:\bbBT_n^{G,\mu}\to\Disp_n^{G,\mu}$ is a gerbe banded by the group scheme $\FD_n^{G,\mu}$ over $\Disp_n^{G,\mu}$.
\eCo

\bproof
Proposition \ref{Pr:bbBT-ttDisp} and Corollary \ref{Co:gerbe-final}.
This proves Theorem \ref{Th:Main}.
\eproof

\bRe
The above reasoning gives a classical construction of a smooth algebraic stack $\ttDisp_n^{G,\mu}$ that is isomorphic to the stack $\bbBT_n^{G,\mu}$, but this does not prove the mod $p$ version of Drinfeld's algebraicity conjecture since the proof of the isomorphism $\bbBT_n^{G,\mu}\cong \ttDisp_n^{G,\mu}$ requires that both sides are known to be classical algebraic stacks. 
\eRe


\end{document}